\documentclass[hidelinks,onefignum,onetabnum]{siamart251216}

\usepackage{amsfonts}
\usepackage{amssymb}
\usepackage{graphicx}
\ifpdf
  \DeclareGraphicsExtensions{.pdf,.png,.jpg,.eps}
\else
  \DeclareGraphicsExtensions{.eps}
\fi

\newcommand{\cT}{\mathcal{T}}

\newcommand{\cV}{\mathcal{V}}

\newcommand{\cZ}{\mathcal{Z}}

\newcommand{\RN}[1]{\uppercase\expandafter{\romannumeral#1}}

\newcommand{\N}{{\mathbb{N}}}

\newcommand{\R}{{\mathbb{R}}}

\newsiamremark{remark}{Remark}
\newsiamremark{assumption}{Assumption}
\crefname{assumption}{Assumption}{Assumptions}

\headers{Burn and density control with prescribed error bounds}{Thomas Berger}

\title{Burn and density control in uncertain tokamak models with
prescribed error bounds\thanks{Funded by the Deutsche Forschungsgemeinschaft (DFG, German Research Foundation) -- Project-ID 544702565.}}

\author{Thomas Berger\thanks{Institut für Mathematik, Martin-Luther-Universität Halle-Wittenberg, 06120~Halle (Saale), Germany (\email{thomas.berger@mathematik.uni-halle.de}).}}

\ifpdf
\hypersetup{
  pdftitle={Burn and density control in uncertain tokamak models with prescribed error bounds},
  pdfauthor={Thomas Berger}
}
\fi

\begin{document}
\tolerance=1200\emergencystretch=1.5em

\maketitle

\begin{abstract}
The stored energy and the electron density of a burning tokamak plasma
must be kept inside prescribed bands, with the auxiliary heating power
and the fuelling rate as actuators.  For the standard zero-dimensional
burn model with a dynamic confinement closure we show that this
two-actuator configuration has vector relative degree $(2,1)$, independently of the
confinement parameters.  We design a model-free funnel controller which
uses an input filter instead of output derivatives and meets all
actuator constraints -- sign, amplitude and rate -- by construction,
the input deficit being returned to the loop as a correction of the
reference.  For a system class containing the model we prove, under the
compatibility assumption that the corrections stay bounded by constants
fixed in advance, that the tracking errors evolve in prescribed funnels
about the corrected references.  Sufficient conditions for the
assumption are given, and the design is illustrated on an ITER
benchmark.
\end{abstract}



%=====================================================================
\section{Introduction}
%=====================================================================

Magnetically confined fusion has entered an engineering phase in which
feedback control is a design constraint rather than an auxiliary
concern.  This paper is about two of its volume-averaged control
objectives, about the structure of the plant that connects them to the
actuators, and about a controller that respects prescribed bounds on
both without knowing that plant.  We first collect a few notions.

\subsection{Plant}\label{sec:primer}

A \emph{tokamak} confines a fully ionised deuterium--tritium plasma in
a toroidal magnetic field \cite{Wes11}.  At temperatures of order
$10\,\mathrm{keV}$ the fusion reaction releases $17.6\,\mathrm{MeV}$
per event, four fifths of it carried out by the neutron and the rest
by the charged $\alpha$ particle, which stays confined and heats the
plasma. A discharge in which this self-heating dominates the external
heating is called \emph{burning}.  On the transport time scale a few
volume-averaged quantities suffice: the densities of fuel, of $\alpha$
\emph{ash} and of neutral fuel, and the stored energy density $E$.
Energy leaves by radiation and by turbulent transport at the rate
$E/\tau_E$, where the \emph{energy confinement time} $\tau_E$ is the
empirical heart of the model, taken from scaling laws regressed over
many machines \cite{IPB99}; neither $\tau_E$ nor the neutral density is
measured in real time.  Two actuators are available and both are
one-sided: auxiliary heating $P_{\mathrm{aux}}\ge0$ and fuelling
$S\ge0$.  The plasma must stay inside operational limits, the penalty
for leaving one being a \emph{disruption} \cite{Hen07}: the normalised
pressure $\beta$ is bounded above by a magnetohydrodynamic stability
limit \cite{Tro84}, the density above by the Greenwald limit
\cite{Gre02} and below by the divertor detachment requirement
\cite{GraSch24}.  In the high-confinement regime the edge relaxes
quasi-periodically in \emph{edge-localised modes} (ELMs) \cite{Zoh96}
and the core in \emph{sawtooth crashes}, each expelling, if unmitigated, a few per cent
of the stored energy in well under a millisecond -- against $\tau_E$ of
order seconds these are pulses, and they are the dominant disturbance
of Section~\ref{sec:dist}.  The controlled phase is the quasi-stationary
\emph{flat-top}, about $400\,\mathrm s$ in the ITER reference scenario
\cite{Hum15}.  What is measured quickly and reliably is what is to be
controlled \cite{Hut02}: the line-averaged electron density $n_e$ by
interferometry and $E$ by a diamagnetic loop.

\subsection{Control objectives}\label{sec:objectives}

The controlled variables are dictated by what must be kept inside a
band and by what can be measured: the stored energy $E$, bounded above
by the $\beta$-limit \cite{Tro84} and pushed up by the burn condition,
with natural actuator $P_{\mathrm{aux}}$; and the electron density
$n_e$, bounded above by the Greenwald limit \cite{Gre02} and below by
the detachment requirement \cite{GraSch24}, with natural actuator $S$.
Both objectives are two-sided and safety critical, which is the
situation funnel control and prescribed performance control (see Section~\ref{sec:modelfree}) are designed for, and the
actuator--output pairing is fixed by the physics.

\subsection{Why model-free}\label{sec:modelfree}

ITER must reach high control performance within roughly $1200$
commissioning discharges, which leaves little time for tuning algorithms
on the device; consequently all ITER control algorithms are required to
be based on control-level models, and \emph{control mathematics} is
explicitly identified as an open research area \cite{Hum15}.  The models
carry irreducible uncertainty where it matters most: the energy
confinement time is taken from empirical scaling laws regressed on
experiments \emph{without} $\alpha$ particle heating, which is why adaptive
estimation was introduced into model-based burn control in the first
place \cite{GraSch22}.  A control law that needs only structural
information -- relative degree, a nonvanishing input gain,
stability of the internal dynamics -- and that guarantees prescribed
error bounds pointwise in time is thus matched to the requirements.
Funnel control, introduced in~\cite{IlcRyaSan02} and further developed in~\cite{BerIlcRya21,BerLeReis18}, is
such a control law. A closely related approach is prescribed performance control (PPC), see \cite{Bechlioulis2008,BechRovi14}. Both aim to guarantee prescribed transient behaviour by constraining the tracking error within suitably defined bounds. The two approaches, however, differ in their controller design, see the survey~\cite{BerIlcRya25} for details.

Learning-based controllers, proposed for this problem \cite{VitMar98}
and experimentally successful on tokamaks \cite{Deg22}, lack a
priori guarantees; a funnel controller can safeguard them by
taking over as a funnel boundary is approached
\cite{LanDen25,LanEtal24,LanKoe25}.

\subsection{Related work}\label{sec:related}

The zero-dimensional balance model is the standard synthesis
model for burn control and has been developed considerably since~\cite{SchKrsTyn03}: nonlinear control and online optimisation using
heating, isotopic fuelling and impurity injection \cite{BoySch14};
robustification against uncertain fuel-line concentrations
\cite{PajSch19}; two-temperature models with adaptive estimation and
optimal allocation of ITER's actuators \cite{GraSch22}; and coupling of
the core to the scrape-off layer and divertor \cite{GraSch23,GraSch24}.
The adequacy of the zero-dimensional description has recently been
examined against profile-resolved closed-loop simulations
\cite{GraSch25}.  Experimentally, feedback regulation of the burn rate
through the deuterium--tritium mixture was demonstrated on JET during
the 2023 D--T campaign \cite{Len25}. Multi-actuator formulations are thus well established in the
literature. The controllers designed in
\cite{BoySch14,GraSch22,GraSch24,PajSch19,SchKrsTyn03} are model-based
state feedback laws.  They assume the full state to be available for
feedback, which for the model of Section~\ref{sec:model} includes the
neutral density and the confinement time, neither of them measured in
real time. The respective stability analysis is carried out without external disturbances, uncertainty being admitted only in parameters of the
plasma model and of the actuators. What is missing is an output-feedback design with prescribed
transient bounds, which is developed in the present paper.

Every ingredient of the desired design exists in the funnel control and PPC literature:
the vector relative degree design \cite{BerLeReis20}, the input filter
\cite{DenSchBer25}, the reference correction \cite{FotRov24} and the
nested-saturation actuator \cite{TraBec23}.  None of them has been
combined with any other before, not even pairwise.  Funnel control under
input constraints has been developed in \cite{Ber24,BerImp26}, where the
funnel boundary is widened while the input saturates; here the funnels
stay prescribed and the reference is corrected instead. Furthermore, we extend the individual components
to fit the present framework: The filter of \cite{DenSchBer25} is
formulated for a constant input gain and is used here with a
state-dependent one; the correction of \cite{FotRov24} is built for a
chain of measured states and is used here on a filter state; and the
filter is a lag that a saturated input winds up, so that the two
mechanisms interact where each was designed to stand alone.  Whether
the assembly is feasible (i.e., it does not violate the funnels) is therefore not clear in advance and
requires a rigorous proof, which is the core part of the present paper.

\subsection{Contributions}

\begin{enumerate}
\item[(I)] We show that the two-actuator configuration has vector
  relative degree $(2,1)$ and a positive diagonal decoupling matrix
  \emph{for structural reasons} (Proposition~\ref{prop:struct}).
\item[(II)] We propose a funnel controller which combines the vector
  relative degree design of \cite{BerLeReis20} with the input filter of
  \cite{DenSchBer25}, so that no derivative of a measured output is
  required, and we prove that the tracking errors evolve in prescribed
  funnels for a system class containing the model
  (Theorem~\ref{thm:main}), which also includes the disturbances typical
  of tokamak models (Section~\ref{sec:dist}).
\item[(III)] Every constraint the ITER actuators impose is met by
  construction rather than verified a posteriori: the sign constraints
  and the installed heating power by returning the input deficit to the
  loop as a correction of the reference, after \cite{FotRov24}, and the
  ramp rate of the heating power by the nested-saturation actuator of
  \cite{TraBec23}, which delivers the amplitude limit with it (Lemma~\ref{lem:act}).
\item[(IV)] We show that the input filter, which is what makes the
  design derivative-free, is a first-order lag driven by the
  \emph{demanded} input and therefore winds up whenever that demand is
  not applied, and that the reference correction is exactly the
  anti-windup it needs: the filter state differs from the one a filter
  driven by the applied input would have exactly by a scaled copy of the
  filtered input deficit, so that it is bounded a priori up to that
  deficit (Lemma~\ref{lem:shadow}).
\item[(V)] The invariance of the funnels does not depend on the sign
of the input gain of the density channel:  Definition~\ref{def:class}
requires that gain only to be nonzero, and the proof of
Theorem~\ref{thm:main} uses only an upper bound of it on the output
box.  The invariance is conditional on bounded
reference corrections and establishes tracking authority for neither
sign (Remark~\ref{rem:cost}).
\end{enumerate}

\subsection{Nomenclature}

For an interval $I\subseteq\mathbb{R}$, $L^1(I;\mathbb{R}^p)$ and
$L^\infty(I;\mathbb{R}^p)$ denote the Lebesgue spaces of measurable
functions $I\to\mathbb{R}^p$ which are integrable and essentially
bounded, respectively, with $\|\cdot\|_\infty$ the essential
supremum; $W^{k,\infty}(I;\mathbb{R}^p)$ is the Sobolev space of
functions with weak derivatives up to order $k$ in
$L^\infty(I;\mathbb{R}^p)$.  $\mathcal{C}^k(I;\mathbb{R}^n)$ is the space
of $k$-times continuously differentiable functions
$I\to\mathbb{R}^n$, $k\in\N_0$; $\mathcal{C}(I;\mathbb{R}^n):=\mathcal{C}^0(I;\mathbb{R}^n)$.  For $r\in\mathbb{R}$ we write
$r_+:=\max\{r,0\}$ and $r_-:=\max\{-r,0\}$.  Densities are in
$\mathrm{m^{-3}}$, energies in $\mathrm{J\,m^{-3}}$, powers in
$\mathrm{W\,m^{-3}}$, temperatures in $\mathrm{keV}$.

%=====================================================================
\section{Model and structural analysis}\label{sec:model}
%=====================================================================

\subsection{Model equations}

We use the zero-dimensional burn model of \cite[Sec.~II.A]{SchKrsTyn03},
augmented by the disturbance terms $d_2,d_3,d_4$ described in
Section~\ref{sec:dist}:
\begin{subequations}\label{eq:model}
	\begin{align}
		\dot n_\alpha(t) &= -\frac{n_\alpha(t)}{k_\alpha\tau_E(t)} + R(t)
		+ d_2(t), \label{eq:na}\\
		\dot n_{DT}(t) &= -\frac{n_{DT}(t)}{k_{DT}\tau_E(t)} - 2R(t)
		+ \frac{n_n(t)}{k_d\tau_E(t)} + d_3(t), \label{eq:ndt}\\
		\dot n_n(t) &= -\frac{n_n(t)}{k_d\tau_E(t)} + S(t), \label{eq:nn}\\
		\dot E(t) &= -\frac{E(t)}{\tau_E(t)} + Q_\alpha R(t)
		- P_{\mathrm{rad}}(t) + P_{\mathrm{aux}}(t) + d_4(t), \label{eq:E}
	\end{align}
\end{subequations}
where
\begin{equation}\label{eq:algebraic}
	\begin{aligned}
		R(t) &= \Big(\tfrac{n_{DT}(t)}{2}\Big)^{2}
		\langle\sigma v\rangle\big(T(t)\big),
		&\quad
		T(t) &= \tfrac{2}{3}\,\frac{E(t)}{2n_{DT}(t)+3n_\alpha(t)},\\
		n_e(t) &= n_{DT}(t)+2n_\alpha(t),
		&\quad
		P_{\mathrm{rad}}(t) &= A_b\big(n_{DT}(t)+4n_\alpha(t)\big)
		n_e(t)\sqrt{T(t)},\\
		P(t) &= Q_\alpha R(t)-P_{\mathrm{rad}}(t)+P_{\mathrm{aux}}(t),
		&\quad
		P_L(t) &= E(t)/\tau_E(t).
	\end{aligned}
\end{equation}
All appearing quantities are described in  Tables~\ref{tab:vars}
and \ref{tab:pars}. The inputs of system~\eqref{eq:model} are the
auxiliary heating power $P_{\mathrm{aux}}(t)$ and the fuelling rate
$S(t)$; the disturbances $d_2,d_3$ are particle sources and $d_4$ is a
power density, and a fourth disturbance $d_1$ enters the confinement
closure \eqref{eq:closure} below.  Solutions of \eqref{eq:model} are
understood in the sense of Carath\'eodory throughout, the disturbances
being merely measurable.
%The second line of \eqref{eq:algebraic}
%carries the conversion $1\,\mathrm{keV}=1.602\cdot10^{-16}\,\mathrm{J}$
%between the units of $E$ and those of $T$ fixed in the nomenclature; we
%suppress it, as is customary.
Ohmic heating is omitted and absorbed into $d_4$. The Maxwellian-averaged D--T reactivity $\langle\sigma v\rangle$ in~\eqref{eq:algebraic},
a function of the temperature alone, is taken from \cite{BosHal92}. The closed-form parametrisation of it presented there exhibits three properties, which are important in the context of this paper and are hence assumed throughout: $\langle\sigma v\rangle$ is nonnegative and locally
Lipschitz on $(0,\infty)$, and it is bounded,
\begin{equation}\label{eq:svmax}
\sigma_v^{\max}:=\sup_{T>0}\langle\sigma v\rangle(T)<\infty .
\end{equation}

\begin{table}[t]
\footnotesize\setlength{\tabcolsep}{3pt}
\begin{minipage}[t]{0.51\textwidth}
\caption{Time-dependent quantities.  The first four rows are the state
of \eqref{eq:model} and the next two the inputs; $\tau_E$ becomes a
fifth state in \eqref{eq:closure}.}
\label{tab:vars}
\centering
\begin{tabular}{@{}lll@{}}
\hline
symbol & meaning & unit\\
\hline
$n_\alpha(t)$ & $\alpha$ particle (ash) density & $\mathrm{m^{-3}}$\\
$n_{DT}(t)$ & ionised D--T fuel density & $\mathrm{m^{-3}}$\\
$n_n(t)$ & neutral fuel, fuelling channel & $\mathrm{m^{-3}}$\\
$E(t)$ & plasma energy density & $\mathrm{J\,m^{-3}}$\\
$P_{\mathrm{aux}}(t)$ & auxiliary heating power & $\mathrm{W\,m^{-3}}$\\
$S(t)$ & fuelling rate & $\mathrm{m^{-3}s^{-1}}$\\
$T(t)$ & plasma temperature & $\mathrm{keV}$\\
$n_e(t)$ & electron density & $\mathrm{m^{-3}}$\\
$R(t)$ & fusion reaction rate & $\mathrm{m^{-3}s^{-1}}$\\
$P_{\mathrm{rad}}(t)$ & bremsstrahlung loss & $\mathrm{W\,m^{-3}}$\\
$P(t)$ & net heating power & $\mathrm{W\,m^{-3}}$\\
$P_L(t)$ & loss power, \eqref{eq:algebraic} & $\mathrm{W\,m^{-3}}$\\
$\tau_E(t)$ & energy confinement time & $\mathrm{s}$\\
$\Lambda(t)$ & particle imbalance, \eqref{eq:Lambda} & $\mathrm{m^{-3}}$\\
$d_1(t)$ & confinement factor, (D1) & --\\
$d_2(t),d_3(t)$ & particle sources, (D2) & $\mathrm{m^{-3}s^{-1}}$\\
$d_4(t)$ & energy bursts, (D3) & $\mathrm{W\,m^{-3}}$\\
\hline
\end{tabular}
\end{minipage}\hfill
\begin{minipage}[t]{0.457\textwidth}
\caption{Constants.  Values are those
for the ITER reference scenario \cite{Shi07}.  The confinement
exponents are those of the ITPA20-IL scaling \cite{Ver21}, obtained from
the ITER-like subset of the revised multi-machine database and
superseding IPB98$(y,2)$; the transport ratios are those of
\cite{Rei90,SchKrsTyn03}; $c_\tau$ is fixed by the operating point, see
Section~\ref{sec:closure}.}
\label{tab:pars}
\centering
\begin{tabular}{@{}lll@{}}
\hline
symbol & meaning & value\\
\hline
$Q_\alpha$ & $\alpha$ birth energy & $3.52\,\mathrm{MeV}$\\
$A_b$ & bremsstrahlung coeff. & $4.85\cdot10^{-37}$\\
$V$ & plasma volume & $830\,\mathrm{m^3}$\\
$c_\tau$ & scaling constant & calibrated\\
$\alpha_P$ & power-degradation exp. & $0.64$\\
$\alpha_n$ & density exponent & $0.15$\\
$k_\alpha$ & ash ratio $\tau_\alpha/\tau_E$ & $7$\\
$k_{DT}$ & fuel ratio $\tau_{DT}/\tau_E$ & $3$\\
$k_d$ & lag ratio $\tau_d/\tau_E$ & $1$\\
$T_c$ & relaxation time in \eqref{eq:closure} & $1\,\mathrm{s}$\\
\hline
\end{tabular}
\end{minipage}
\end{table}

The ratios $k_\alpha,k_{DT},k_d$ tie the particle time scales to
$\tau_E$; $k_\alpha$ is the figure of merit of helium exhaust, required
to stay below about $10$ \cite{Rei90}, and $k_\alpha=7$, $k_{DT}=3$ are
the values of the modelling tradition.  %The model carries no impurity
%species; an impurity state treated as a disturbance would leave
%Proposition~\ref{prop:struct} intact, whereas impurity injection as a
%third input would add an input direction and require the structural
%analysis to be redone.

\subsection{Dynamic confinement closure}\label{sec:closure}

The empirical component of the model is the confinement time.  It is
customary to insert the steady-state correlation instantaneously,
$\tau_E(t)=c_\tau n_e(t)^{\alpha_n}P(t)^{-\alpha_P}$.  We do not do
this, for a physical reason: a correlation regressed on stationary
discharges cannot be expected to hold instantaneously, and the
confinement cannot respond faster than the transport time scale.  We
therefore close the model dynamically, and with the \emph{loss} power
$P_L=E/\tau_E$, against which the scalings are regressed and
by which $\tau_E$ is defined \cite{IPB99,Ver21}, rather than the net
heating power $P$; by \eqref{eq:E}, $P_L=P+d_4-\dot E$, which is
$P-\dot E$ in the absence of the disturbance $d_4$ and coincides with
$P$ only in stationary discharges without $d_4$:
\begin{equation}\label{eq:closure}
\dot\tau_E(t) = \frac{1}{T_c}\Big(
  d_1(t)\,c_\tau\,n_e(t)^{\alpha_n}P_L(t)^{-\alpha_P} - \tau_E(t)\Big),
\end{equation}
with a relaxation time $T_c>0$, so that $\tau_E$ becomes a fifth state,
and with the multiplicative disturbance $d_1$ of Section~\ref{sec:dist},
which is $1$ for the nominal scaling.
No experimental value of $T_c$ is available; it is bounded below by
the saturation time of the turbulence and above by the profile
relaxation time, of order $\tau_E$ itself, and we use $T_c=1\,\mathrm s$.
At any equilibrium with $d_1\equiv1$ and $d_4\equiv0$ we have $P_L=P$, so
\eqref{eq:closure} has exactly the equilibria of the algebraic closure:
it is a modelling
refinement that shares its equilibria with the standard closure, not a
competing model.

We fix the operating point of system \eqref{eq:model}, \eqref{eq:closure},
with $d_1\equiv1$ and $d_2\equiv d_3\equiv d_4\equiv0$, once and for
all as the ITER reference
scenario \cite{Shi07}: a line-averaged density
$\bar n_e=1.01\cdot10^{20}\,\mathrm{m^{-3}}$, i.e.\ $85\%$ of the
Greenwald limit, a temperature $\bar T=8.4\,\mathrm{keV}$ and an ash
fraction $\bar n_\alpha/\bar n_e=3.7\%$, whence
$\bar E=4.00\cdot10^{5}\,\mathrm{J\,m^{-3}}$.  Calibrating $c_\tau$ so
that \eqref{eq:model}, \eqref{eq:closure} has this equilibrium yields
$\bar\tau=3.48\,\mathrm{s}$ and $\bar P_{\mathrm{aux}}=36.4\,$MW,
against the $3.7\,\mathrm{s}$ that \cite{Shi07} predicts for that
scenario and the $50\,$MW to which stationary operation at $Q=10$
confines the auxiliary power \cite{Sip15}.  This checks the calibration at the operating point, but
validates neither the transients of \eqref{eq:closure} nor the value of $T_c$; the analysis of
Sections~\ref{sec:model}--\ref{sec:apriori} holds for every $T_c>0$.
The equilibrium is held by the nominal inputs
$S\equiv\bar S=9.26\cdot10^{18}\,\mathrm{m^{-3}\,s^{-1}}$ and
$P_{\mathrm{aux}}\equiv\bar P:=\bar P_{\mathrm{aux}}/V
=4.39\cdot10^{4}\,\mathrm{W\,m^{-3}}$, the latter a power density, the
unit in which $P_{\mathrm{aux}}$ enters \eqref{eq:E}.

\subsection{Disturbances}\label{sec:dist}

A burning plasma is subject to a
variety of disturbances.  The deviation between the nominal model and a
discharge is carried by the four terms $d_1,\dots,d_4$ in
\eqref{eq:model}, \eqref{eq:closure}, which model three mechanisms.

\noindent
\textbf{(D1) Confinement uncertainty.}  The correlation in
\eqref{eq:closure} has a log-linear uncertainty band of $15$--$20\%$
\cite{IPB99,Ver21} and a prefactor that drops at the L--H threshold
\cite{SchKrsTyn03}.  Both are modelled by the factor $d_1$ in
\eqref{eq:closure}, with $d_1\in L^{\infty}(\mathbb{R}_{\ge0};\mathbb{R})$
such that $\operatorname{ess\,inf}_{t\ge 0} d_1(t)>0$.

\noindent
\textbf{(D2) Particle sources and sinks.}  Wall recycling, neutral
influx and isotopic imbalance enter \eqref{eq:na}, \eqref{eq:ndt}
additively as the source terms $d_2,d_3\in
W^{1,\infty}(\mathbb{R}_{\ge0};\mathbb{R})$, with $d_2,d_3\ge0$ so that no
density is drawn negative.  These two are the only disturbances of
which more than boundedness is asked; the reason is that the density
channel is differentiated twice, see Proposition~\ref{prop:struct}
below.

\noindent
\textbf{(D3) Energy bursts.}  Edge-localised modes and sawtooth crashes
expel a fraction of the stored energy on a time scale far below
$\tau_E$. Together with impurity radiation and the neglected ohmic
heating they form the term $d_4$, which enters the energy balance
\eqref{eq:E} additively as a power density.  We only require
$d_4\in L^\infty(\mathbb{R}_{\ge0};\mathbb{R})$: a
pulse train is discontinuous and not of bounded
variation, but bounded and measurable.

\subsection{The two-actuator configuration}\label{sec:mimo}

We now take
\begin{equation}\label{eq:io}
u=(u_1,u_2)=(S,P_{\mathrm{aux}}),\qquad
y=(y_1,y_2)=(n_e,E),
\end{equation}
as input and output vector, the pairing being the one fixed in
Section~\ref{sec:objectives}.  Write
$x:=(n_\alpha,n_{DT},n_n,E,\tau_E)$, $d:=(d_1,d_2,d_3,d_4)$ and
$\dot x(t)=f(x(t),u(t),d(t))$ for \eqref{eq:model}, \eqref{eq:closure}, where $f:\mathcal{X}\times \R^2\times \R^4 \to \R^5$, with
\begin{equation}\label{eq:Xdef}
\mathcal{X}:=\{n_\alpha\ge0,\ n_{DT}\ge0,\ n_\alpha+n_{DT}>0,\
 n_n\ge0,\ E>0,\ \tau_E>0\} \subseteq\R^5,
\end{equation}
such that $f(\cdot,u,d)$ is the restriction of a locally Lipschitz map
defined on an open neighbourhood, for every $u$ and $d$; the condition
$n_\alpha+n_{DT}>0$ is what $T$ in \eqref{eq:algebraic} requires and
$E>0$ is required to make $P_L^{-\alpha_P}$ in~\eqref{eq:closure}
well-defined for $P_L = E/\tau_E$ as in \eqref{eq:algebraic}.  Finally,
$y(t)=h(x(t))$ with
$h:\mathcal{X}\to\R^2,\ x\mapsto (n_{DT}+2n_\alpha,\,E)$.  Put
\[
\mathcal{G}_i:\mathcal{X}\times \R^2\times\R^4\to \R,\
 (x,u,d) \mapsto h_i'(x) f(x,u,d),\quad i=1,2,
\]
where $h_i'$ is the Jacobian of $h_i$, so that
$\dot y_i(t)=\mathcal{G}_i(x(t),u(t),d(t))$ for almost all $t$ along
every solution.  Two exact identities determine the structure.

\begin{proposition}\label{prop:struct}
For all $(x,u,d)\in\mathcal{X}\times\mathbb{R}^2\times\mathbb{R}^4$,
\begin{equation}\label{eq:Lambda}
\mathcal{G}_1(x,u,d)=\frac{\Lambda(x)}{\tau_E}+2d_2+d_3,\qquad
\Lambda(x):=\frac{n_n}{k_d}-\frac{n_{DT}}{k_{DT}}
        -\frac{2n_\alpha}{k_\alpha};
\end{equation}
in particular $\mathcal{G}_1$ does not depend on $u$, and with
$\mathcal{G}_1^0(x):=\Lambda(x)/\tau_E$ and
$\mathcal{G}_{11}(x,u,d):=\frac{\partial\mathcal{G}_1^0}{\partial x}(x)\,
f(x,u,d)$ every solution with $d_2,d_3$ absolutely continuous satisfies
\begin{equation}\label{eq:ddoty1}
\ddot y_1(t)=\mathcal{G}_{11}\big(x(t),u(t),d(t)\big)
 +2\dot d_2(t)+\dot d_3(t)
\end{equation}
for almost all $t$.  Moreover,
\begin{equation}\label{eq:Gamma}
\frac{\partial(\mathcal{G}_{11},\mathcal{G}_2)}{\partial u}(x,u,d)
 =\hat\Gamma(x):=\operatorname{diag}
 \Big(\tfrac{1}{k_d\tau_E},\,1\Big)
\end{equation}
for all $(x,u,d)$, which is invertible and does not depend on $d$.
Since $f$ is affine in $u$, so are $\mathcal{G}_{11}$ and $\mathcal{G}_2$:
\begin{equation}\label{eq:affine}
\mathcal{G}_{11}(x,u,d)=\mathcal{G}_{11}(x,0,d)+\frac{u_1}{k_d\tau_E},
\quad
\mathcal{G}_2(x,u,d)=\mathcal{G}_2(x,0,d)+u_2 .
\end{equation}
\end{proposition}

\begin{proof}
$h_1'=(2,1,0,0,0)$, so $\mathcal{G}_1$ is twice the first plus the
second component of~$f$. Adding \eqref{eq:ndt} and twice \eqref{eq:na}
gives \eqref{eq:Lambda}.  Since $u_1=S$ enters $f$ only through
\eqref{eq:nn} and $u_2=P_{\mathrm{aux}}$ only through \eqref{eq:E},
and the right-hand side of \eqref{eq:closure} depends on $x$ and $d_1$
alone, $\mathcal{G}_1$ does not depend on $u$.  Along a solution,
$y_1=h_1(x)$ is absolutely continuous with
$\dot y_1=\mathcal{G}_1^0(x)+2d_2+d_3$ a.e.; $\mathcal{G}_1^0$ is
$\mathcal{C}^1$ on $\mathcal{X}$ and $x$ is absolutely continuous, so
$\mathcal{G}_1^0(x(\cdot))$ is absolutely continuous with derivative
$\mathcal{G}_{11}(x(\cdot),u(\cdot),d(\cdot))$ a.e., and \eqref{eq:ddoty1} follows when $d_2,d_3$
are absolutely continuous.  Since
$\frac{\partial \mathcal{G}_1^0}{\partial x}(x)=(-2/(k_\alpha\tau_E),-1/(k_{DT}\tau_E),
1/(k_d\tau_E),0,-\Lambda/\tau_E^2)$, it follows that
$\frac{\partial\mathcal{G}_{11}}{\partial u_1}(x,u,d)=1/(k_d\tau_E)>0$ and
$\frac{\partial\mathcal{G}_{11}}{\partial u_2}(x,u,d)=0$.  Finally
$h_2'=(0,0,0,1,0)$ gives $\mathcal{G}_2=f_4$, whence
$\frac{\partial\mathcal{G}_2}{\partial u_1}(x,u,d)=0$ and
$\frac{\partial\mathcal{G}_2}{\partial u_2}(x,u,d)=1$.  Equation
\eqref{eq:affine} follows because $f$ is affine in $u$ and
$\frac{\partial \mathcal{G}_1^0}{\partial x}$ does not depend on $u$.
\end{proof}

Observe that the decoupling matrix $\hat \Gamma(x)$ defined in~\eqref{eq:Gamma} is globally diagonal, so the
\emph{direct} influence of each actuator is confined to one output.  This does not mean that both channels are independent, the
drifts $\mathcal{G}_{11}$ and $\mathcal{G}_2$ induce a coupling of the outputs. As a consequence of the above result, system~\eqref{eq:model}, \eqref{eq:closure}
 has \emph{vector relative degree} $(2,1)$, a notion that we will make more precise after Definition~\ref{def:class}.

\begin{remark}[Regularity of the disturbances]\label{rem:D2}
Proposition~\ref{prop:struct} exhibits where the different requirements on
$d_1,\dots,d_4$ in Section~\ref{sec:dist} come from.  By
\eqref{eq:Lambda} the particle sources $d_2$ and $d_3$ enter $\dot y_1$ directly, and
\eqref{eq:ddoty1} carries their derivatives, which is why (D2) asks for
$W^{1,\infty}$.  The factor $d_1$ and the bursts $d_4$ enter
$\mathcal{G}_{11}$ and $\mathcal{G}_2$ only through $f$, which is never
differentiated with respect to $t$, so $L^\infty$ suffices.
\end{remark}

\begin{remark}[Invariance of $\mathcal{X}$]\label{rem:physical}
Let the inputs be nonnegative, $S(t)\ge0$ and $P_{\mathrm{aux}}(t)\ge0$,
as the physics requires, and let $d_2,d_3\ge0$ as in (D2).  Then four of
the six conditions in \eqref{eq:Xdef} are positively invariant for
\eqref{eq:model}, \eqref{eq:closure}: the sign conditions on the
densities because $\dot n_n|_{n_n=0}=S\ge0$,
$\dot n_{DT}|_{n_{DT}=0}=n_n/(k_d\tau_E)+d_3\ge0$ and
$\dot n_\alpha|_{n_\alpha=0}=R+d_2\ge0$, and $\tau_E>0$ because the
first term of $\dot\tau_E$ in \eqref{eq:closure} is proportional to
$\tau_E^{\alpha_P}$ with $\alpha_P\in(0,1)$ and dominates $-\tau_E$ as
$\tau_E\downarrow0$.  The remaining two are not necessarily positively invariant: $\dot E|_{E=0}
=P_{\mathrm{aux}}+d_4$ can be negative, since the energy bursts $d_4$
are, and for nonnegative densities $n_\alpha+n_{DT}>0$ holds exactly when the first output
$n_e=n_{DT}+2n_\alpha$ is positive. Keeping the outputs $(n_e,E)$
within a desired range is the task of the controller design.
\end{remark}

%=====================================================================
\section{System class}\label{sec:class}
%=====================================================================

As mentioned above, Proposition~\ref{prop:struct} uncovers that system~\eqref{eq:model}, \eqref{eq:closure} has a vector relative degree, which
places it in the setting of \cite{BerLeReis20}: square, ordered, with $q=m=2$, so that the
algebraic block of \cite[Eq.~(18)]{BerLeReis20} is absent.  The class
below differs from that reference in three respects, each dictated by
the plant: the input gain of the first channel is state dependent and
is carried by an operator, the input gain of the second channel is a
constant, and the internal dynamics is defined only on output histories
that stay in a prescribed set $\mathcal{Z}$, which becomes the
parameter of the class.

\begin{definition}\label{def:class}
Let $\mathcal{Z}\subseteq\mathcal{C}(\mathbb{R}_{\ge0};\mathbb{R}^2)$
be a set of \emph{admissible output histories}.  The
class $\mathcal{B}(\mathcal{Z})$ consists of all systems
\begin{equation}\label{eq:class}
\begin{aligned}
\ddot y_1(t) &= F_1(t) + \gamma(t)\,u_1(t),\\
\dot y_2(t) &= F_2(t) + a_2\,u_2(t),
\end{aligned}
\end{equation}
with a constant $a_2>0$, for which there exist $p,q\in\mathbb{N}$, a
disturbance $d\in L^\infty(\mathbb{R}_{\ge0};\mathbb{R}^p)$, an operator
\[
\mathcal{T}:\mathcal{Z}\to
 L^\infty_{\mathrm{loc}}(\mathbb{R}_{\ge0};\mathbb{R}^q),
\]
which is fixed together with $d$ and may depend on it, an open set
$\mathcal{V}\subseteq\mathbb{R}^2\times\mathbb{R}^p\times\mathbb{R}^q$,
from which we form the open sets
\begin{equation}\label{eq:domains}
\begin{aligned}
\mathcal{V}_1&:=\big\{(y_1,\vartheta,y_2,\delta,w)\ \big|\
 \vartheta\in\mathbb{R},\ (y_1,y_2,\delta,w)\in\mathcal{V}\big\},\\
\mathcal{W}&:=\big\{w\ \big|\ (y_1,y_2,\delta,w)\in\mathcal{V}
 \text{ for some }(y_1,y_2,\delta)\big\},
\end{aligned}
\end{equation}
and functions
$f_1:\mathcal{V}_1\to\mathbb{R}$,
$f_2:\mathcal{V}\to\mathbb{R}$ and
$g:\mathcal{W}\to\mathbb{R}$ such that
\begin{equation}\label{eq:Fdef}
\begin{aligned}
F_1(t)&=f_1\big(y_1(t),\dot y_1(t),y_2(t),d(t),
 \mathcal{T}(\zeta)(t)\big),\\
F_2(t)&=f_2\big(y_1(t),y_2(t),d(t),\mathcal{T}(\zeta)(t)\big),\\
\gamma(t)&=g\big(\mathcal{T}(\zeta)(t)\big),
\end{aligned}
\end{equation}
and
\begin{enumerate}
\item[(B1)] $f_2$ and $g$ are continuous, and $g\ne0$ on
  $\mathcal{W}$;
\item[(B2)] $f_1$ is continuous and, for some continuous
  $c:\mathcal{V}\to\mathbb{R}_{\ge0}$ and some constant $\bar B\ge0$,
  \[
  |f_1(y_1,\hat y_1,y_2,\delta,w)|\le c(y_1,y_2,\delta,w)+\bar B\,|\hat y_1|
  \]
  for all $\hat y_1\in\mathbb{R}$ and all
  $(y_1,y_2,\delta,w)\in\mathcal{R}$, where
  \[
  \mathcal{R}:=\bigcup_{\zeta\in\mathcal{Z}}
   \operatorname{ess\,im}\big(\zeta,d,\mathcal{T}(\zeta)\big)
   \subseteq\mathcal{V}
  \]
  is the set of values attained along admissible histories: for
  $\zeta\in\mathcal{Z}$, $(\zeta,d,\mathcal{T}(\zeta))$ denotes the
  measurable function $t\mapsto(\zeta(t),d(t),\mathcal{T}(\zeta)(t))$ on
  $\mathbb{R}_{\ge0}$, and the essential image $\operatorname{ess\,im}\,g$
  of a measurable function $g$ is the set of points every neighbourhood
  of which has a preimage of positive measure, so that $g(t)\in
  \operatorname{ess\,im}\,g$ for almost all $t$;
\item[(B3)] $\mathcal{T}$ is \emph{causal},
  \[
  \zeta|_{[0,t]}=\tilde{\zeta}|_{[0,t]}\ \Longrightarrow\
   \mathcal{T}(\zeta)|_{[0,t]}=\mathcal{T}(\tilde{\zeta})|_{[0,t]}
   \quad\text{for all }t\ge0,\ \zeta,\tilde{\zeta}\in\mathcal{Z},
  \]
  \emph{locally Lipschitz}, i.e.\ for every $t\ge0$ and every
  $\tilde{\zeta}\in\mathcal{C}([0,t];\mathbb{R}^2)$ there are
  $c_0,\varrho,\varsigma>0$ such that all
  $\zeta_1,\zeta_2\in\mathcal{Z}$ with $\zeta_j|_{[0,t]}=\tilde{\zeta}$ and
  $\|\zeta_j(s)-\tilde{\zeta}(t)\|<\varrho$ on $[t,t+\varsigma]$ satisfy
  \[
  \operatorname*{ess\,sup}_{s\in[t,t+\varsigma]}
   \|\mathcal{T}(\zeta_1)(s)-\mathcal{T}(\zeta_2)(s)\|
   \le c_0\sup_{s\in[t,t+\varsigma]}\|\zeta_1(s)-\zeta_2(s)\| ,
  \]
  and \emph{bounded-input bounded-output}: for every $c_1>0$ there is
  a compact set $\mathcal{K}\subseteq\mathcal{V}$ such that every
  $\zeta\in\mathcal{Z}$ with $\|\zeta\|_\infty\le c_1$ satisfies
  \[
  \big(\zeta(t),d(t),\mathcal{T}(\zeta)(t)\big)\in\mathcal{K}
  \quad\text{for almost all $t\ge0$.}
  \]
\end{enumerate}
\end{definition}

Only $\mathcal{Z}$ parametrises $\mathcal{B}$; the dimensions, the
constant $a_2$, the disturbance, the operator, the domain, the three
functions and the growth constant $\bar B$ of (B2) belong to the
individual system and are free within the class.  The set
$\mathcal{Z}$ is the interface between plant and design: it is the
set of output histories on which the plant's internal dynamics is
known to be well behaved, and the control design of Section~\ref{sec:main} will
have to keep the outputs inside it.  The domain is what lets the
class carry functions with poles, as the
plant of Section~\ref{sec:model} has at $\tau_E=0$; (B3) keeps the
signals away from them.

Note that, by causality,
$\mathcal{T}(\zeta)|_{[0,t]}$ depends on $\zeta|_{[0,t]}$ alone; a
history $\zeta\in\mathcal{C}([0,\omega);\mathbb{R}^2)$ whose
restrictions $\zeta|_{[0,t]}$, $t<\omega$, all admit an extension in
$\mathcal{Z}$ therefore has a well-defined image
$\mathcal{T}(\zeta)\in L^\infty_{\mathrm{loc}}([0,\omega);\mathbb{R}^q)$,
independent of the extensions chosen. In this way the operator acts on
solutions on $[0,\omega)$ throughout.

The class \eqref{eq:class} is tailored so that it encompasses~\eqref{eq:model}, \eqref{eq:closure}
for all possible parameter configurations and disturbances~-- but otherwise as simple as possible. Therefore, the second input gain is a
constant, the off-diagonal gains vanish, and $\gamma$ is a function of
the internal variables alone.  No sign of $\gamma$ is required
(Remark~\ref{rem:cost}).

The system \eqref{eq:class} is said to have
\emph{vector relative degree} $(2,1)$ with \emph{decoupling matrix}
$\operatorname{diag}(\gamma,a_2)$, which is nonsingular by (B1).  The
operator $\mathcal{T}$ models
the \emph{internal dynamics} of the system and (B3) asks of that part
no more than bounded-input bounded-output stability. The operator~$\mathcal{T}$ is
 a standard modelling tool in funnel control (also capable of describing infinite-dimensional internal dynamics),
 see~\cite{BerIlcRya25}.

\subsection{The plant as a member of the system class}

Placing the plant~\eqref{eq:model}, \eqref{eq:closure} in the class \eqref{eq:class} takes four objects -- a box of output
histories, an extension of the model's nonlinear terms, a change of
coordinates and a handful of constants.  We define them first and state
membership afterwards.  The disturbances $d_1,\dots,d_4$ are fixed
throughout as in Section~\ref{sec:dist}.

\emph{The box.}  Fix numbers $0<\underline y_i\le\overline y_i$,
$i=1,2$, and let
\[
\mathcal{Z}:=\big\{\zeta=(y_1,y_2)\in
 \mathcal{C}(\mathbb{R}_{\ge0};\mathbb{R}^2)\ \big|\
 \underline y_i\le y_i(t)\le\overline y_i\ \text{for all }t\ge0\big\} .
\]

\emph{The clamped extension.}  The reaction rate $R$ and the radiation
$P_{\mathrm{rad}}$ in \eqref{eq:algebraic} are defined only where $T>0$.
To extend them to all values of the ash density, let
$\mathbb{D}:=(0,\infty)\times(0,\infty)\times\mathbb{R}$,
whose points $\mathbf{y}=(y_1,y_2,a)$ stand for the values of $n_e=h_1(x)$,
$E=h_2(x)$ and $n_\alpha$, where $h$ is defined in Section~\ref{sec:mimo}, and define, for $\mathbf{y}=(y_1,y_2,a)\in\mathbb{D}$,
\begin{equation}\label{eq:Rext}
\begin{aligned}
a_*(\mathbf{y})&:=\min\{\max\{a,0\},\,y_1/2\},\qquad
n_{DT}^*(\mathbf{y}):=y_1-2a_*(\mathbf{y}),\\
T_*(\mathbf{y})&:=\frac{2}{3}\,\frac{y_2}{2y_1-a_*(\mathbf{y})},\qquad
R_*(\mathbf{y}):=\tfrac14\,n_{DT}^*(\mathbf{y})^2\,\langle\sigma v\rangle\big(T_*(\mathbf{y})\big),
\\
P_{\mathrm{rad}}^*(\mathbf{y})&:=A_b\,\big(n_{DT}^*(\mathbf{y})+4a_*(\mathbf{y})\big)\,y_1
 \sqrt{T_*(\mathbf{y})} .
\end{aligned}
\end{equation}
Since $0\le a_*(\mathbf{y})\le y_1/2$, we have $n_{DT}^*(\mathbf{y})\in[0,y_1]$ and
$2y_1-a_*(\mathbf{y})\ge\tfrac32y_1>0$, so $T_*(\mathbf{y})>0$; as
$\langle\sigma v\rangle$ is nonnegative on $(0,\infty)$, \eqref{eq:Rext} defines maps
\[
a_*,\ n_{DT}^*:\mathbb{D}\to[0,\infty),\qquad
T_*:\mathbb{D}\to(0,\infty),\qquad
R_*,\ P_{\mathrm{rad}}^*:\mathbb{D}\to[0,\infty).
\]
All five are locally Lipschitz: $a_*$ is globally Lipschitz as a
composition of $\min$, $\max$ and linear maps, hence so is
$n_{DT}^*$; $T_*$ is a quotient of locally Lipschitz maps whose
denominator is locally bounded away from zero; and $R_*$,
$P_{\mathrm{rad}}^*$ are products and compositions of these with
$\langle\sigma v\rangle$ and $\sqrt{\cdot}$, both locally Lipschitz on
$(0,\infty)$.  For $x=(n_\alpha,n_{DT},n_n,E,\tau_E)\in\mathcal{X}$ the
point $\mathbf{y}_x:=(h_1(x),h_2(x),n_\alpha)$ lies in $\mathbb{D}$ and satisfies
$0\le n_\alpha\le h_1(x)/2$, hence $a_*(\mathbf{y}_x)=n_\alpha$,
$n_{DT}^*(\mathbf{y}_x)=n_{DT}$ and $T_*(\mathbf{y}_x)=T$, and $R_*(\mathbf{y}_x)$,
$P_{\mathrm{rad}}^*(\mathbf{y}_x)$ are exactly $R$ and $P_{\mathrm{rad}}$ of
\eqref{eq:algebraic} at $x$.
On the open set
\[
\mathcal{O}:=\{x\in\mathbb{R}^5\mid\tau_E>0,\ E>0,\ n_{DT}+2n_\alpha>0\}
 \supset\mathcal{X}
\]
define $f^*:\mathcal{O}\times\mathbb{R}^2\times\mathbb{R}^4\to\mathbb{R}^5$,
$(x,u,d)\mapsto f^*(x,u,d)$, by the formulas of the right-hand side
$f(x,u,d)$ of \eqref{eq:model}, \eqref{eq:closure}, with $R$ and
$P_{\mathrm{rad}}$ replaced by $R_*(\mathbf{y}_x)$ and $P_{\mathrm{rad}}^*(\mathbf{y}_x)$,
$\mathbf{y}_x=(h_1(x),h_2(x),n_\alpha)\in\mathbb{D}$, which is well defined for
every $x\in\mathcal{O}$.  Then $f^*$ is locally Lipschitz -- the conditions
$E>0$ and $n_{DT}+2n_\alpha>0$ in $\mathcal{O}$ are what $T_*$ and
$n_e^{\alpha_n}$ need -- and coincides with $f$ on
$\mathcal{X}\times\mathbb{R}^2\times\mathbb{R}^4$.  With
$\mathcal{G}_1^0(x)=\Lambda(x)/\tau_E$ from Proposition~\ref{prop:struct}
put, for $(x,u,d)\in\mathcal{O}\times\mathbb{R}^2\times\mathbb{R}^4$,
\[
\begin{aligned}
\mathcal{G}^*_2(x,u,d)&:=\nabla h_2(x)\,f^*(x,u,d),\\
\mathcal{G}^*_{11}(x,u,d)&:=\nabla\mathcal{G}_1^0(x)\,f^*(x,u,d),
\end{aligned}
\]
which coincide with $\mathcal{G}_2,\mathcal{G}_{11}$ of
Proposition~\ref{prop:struct} on $\mathcal{X}\times\mathbb{R}^2\times
\mathbb{R}^4$.

\emph{The coordinates.}  For $x=(n_\alpha,n_{DT},n_n,E,\tau_E)\in\mathbb{R}^4\times (0,\infty)$ let
$\eta(x):=(n_\alpha,\tau_E)$ and
$\Xi(x):=\big(h_1(x),\mathcal{G}_1^0(x),h_2(x),\eta(x)\big)$,
and for $(y_1,z,y_2,n_\alpha,\tau_E)\in\mathbb{R}^4\times(0,\infty)$ define the smooth map
\begin{equation}\label{eq:nninv}
  \iota(y_1,z,y_2,n_\alpha,\tau_E)
  :=\left(n_\alpha,\,y_1-2n_\alpha,\,k_d\Big(\tau_E z+\frac{y_1-2n_\alpha}{k_{DT}}
  +\frac{2n_\alpha}{k_\alpha}\Big),\,y_2,\,\tau_E\right).
\end{equation}

\emph{The constants.}  Put
$\underline d_1:=\operatorname*{ess\,inf}_{t\ge0}d_1(t)>0$ and
$\overline d_1:=\operatorname*{ess\,sup}_{t\ge0}d_1(t)<\infty$, so that
$d_1(t)\in[\underline d_1,\overline d_1]$ for almost all $t$, and, for
numbers $\delta,y_1,y_2>0$,
$\mathcal{Q}(\delta,y_1,y_2):=\delta\,c_\tau\,y_1^{\alpha_n}y_2^{-\alpha_P}$,
so that by $P_L=E/\tau_E$ the right-hand side of \eqref{eq:closure} is
$\big(\mathcal{Q}(d_1(t),n_e(t),E(t))\,\tau_E(t)^{\alpha_P}-\tau_E(t)\big)
/T_c$.  Since $\mathcal{Q}$ is increasing in $\delta$ and $y_1$ and
decreasing in $y_2$,
\[
\underline{\mathcal Q}:=\mathcal{Q}(\underline d_1,\underline y_1,
 \overline y_2)\le\mathcal{Q}(\delta,y_1,y_2)\le
 \mathcal{Q}(\overline d_1,\overline y_1,\underline y_2)
 =:\overline{\mathcal Q}
\]
for all $\delta\in[\underline d_1,\overline d_1]$ and
$y_i\in[\underline y_i,\overline y_i]$.  Fix an initial value
$\eta^0=(n_\alpha^0,\tau_E^0)\in[0,\infty)\times(0,\infty)$ of the
internal state and define, with $\sigma_v^{\max}$ from \eqref{eq:svmax},
\[
\tau_{\min}:=\min\{\tau_E^0,\underline{\mathcal Q}^{\frac{1}{1-\alpha_P}}\},
\quad
\tau_{\max}:=\max\{\tau_E^0,\overline{\mathcal Q}^{\frac{1}{1-\alpha_P}}\},
\quad
\Sigma:=\tfrac14\overline y_1^2\sigma_v^{\max}+\|d_2\|_\infty ,
\]
and
\begin{equation}\label{eq:Bbar}
\bar{B}:=\frac{1}{\tau_{\min}}\Big(\frac{1}{k_d}+\frac{1}{k_{DT}}\Big)
 +\frac{\Theta}{T_c},\qquad
\Theta:=\max\big\{
 \overline{\mathcal Q}\,\tau_{\min}^{\alpha_P-1}-1,\;
 1-\underline{\mathcal Q}\,\tau_{\max}^{\alpha_P-1}\big\} ,
\end{equation}
both entries of $\Theta$ being nonnegative.

\begin{proposition}\label{prop:membership}
Let $T_c>0$ and let $\mathcal{Z}$, $d_1,\ldots,d_4$ and $\eta^0$ be
as above.
\begin{enumerate}
\item[(i)] \emph{Internal dynamics.}  For every
  $\zeta=(y_1,y_2)\in\mathcal{Z}$ the initial-value problem
  \begin{equation}\label{eq:intdyn}
  \begin{aligned}
  \dot n_\alpha(t) &= -\frac{n_\alpha(t)}{k_\alpha\tau_E(t)}
    + R_*\big(y_1(t),y_2(t),n_\alpha(t)\big) + d_2(t),\\
  \dot\tau_E(t) &= \frac{1}{T_c}\Big(
    \mathcal{Q}\big(d_1(t),y_1(t),y_2(t)\big)\,\tau_E(t)^{\alpha_P}
    -\tau_E(t)\Big),\\
  \big(n_\alpha(0),\tau_E(0)\big)&=\eta^0,
  \end{aligned}
  \end{equation}
  has a unique global solution $\eta(\cdot\,;\zeta)=(n_\alpha,\tau_E)$,
  which satisfies
  \begin{equation}\label{eq:bibobounds}
  \tau_E(t)\in[\tau_{\min},\tau_{\max}],\qquad
  0 \le n_\alpha(t)\le\max\{n_\alpha^0,\,k_\alpha\tau_{\max}\Sigma\}
  \qquad\text{for all $t\ge0$.}
  \end{equation}
  The operator
  $\cT:\cZ\to L^\infty_{\mathrm{loc}}(\mathbb{R}_{\ge0};\mathbb{R}^2)$, $\zeta\mapsto \eta(\cdot\,;\zeta)$,
  satisfies (B3).
\item[(ii)] \emph{Coordinates.}  The map $\Xi:\R^4\times (0,\infty) \to \mathbb{R}^4\times(0,\infty)$ is a diffeomorphism,
  and its inverse is $\iota$.  Let the inputs be nonnegative as in
  Remark~\ref{rem:physical}, so that every solution~$x$ of
  \eqref{eq:model}, \eqref{eq:closure} with $x(0)\in\mathcal{X}$, $\eta(x(0))=\eta^0$ and
  $(h_1(x),h_2(x))\in\mathcal{Z}$ stays in $\mathcal{X}$. Then, written in the coordinates
  $(y_1,z,y_2,\eta)=\Xi(x)$, $\eta=(n_\alpha,\tau_E)$, the plant
  \eqref{eq:model}, \eqref{eq:closure}, \eqref{eq:io} takes the form
  \eqref{eq:class}, \eqref{eq:Fdef} with $p=6$, $q=2$, $a_2=1$,
  $\mathcal{T}$ from (i), the disturbance
  $d:=(d_1,d_2,d_3,d_4,\dot d_2,\dot d_3)
   \in L^\infty(\mathbb{R}_{\ge0};\mathbb{R}^6)$,
  the open set
  \begin{equation}\label{eq:Vplant}
  \mathcal{V}:=\big\{(y_1,y_2,\delta,w)\in\mathbb{R}^2\times\mathbb{R}^6
   \times\mathbb{R}^2\ \big|\ y_1,y_2>0,\ w_2>0\big\},
  \end{equation}
  the sets formed from it as in \eqref{eq:domains},
  \[
  \mathcal{V}_1=\big\{(y_1,\vartheta,y_2,\delta,w)\in\mathbb{R}^3\times
   \mathbb{R}^6\times\mathbb{R}^2\ \big|\ y_1,y_2>0,\ w_2>0\big\},
  \quad
  \mathcal{W}=\mathbb{R}\times(0,\infty),
  \]
  and the functions
  $f_1:\mathcal{V}_1\to\mathbb{R}$, $f_2:\mathcal{V}\to\mathbb{R}$,
  $g:\mathcal{W}\to\mathbb{R}$ given by
  \begin{equation}\label{eq:fdef}
  \begin{aligned}
  f_1(y_1,\vartheta,y_2,\delta,w)
   &:=\mathcal{G}^*_{11}\big(\iota(y_1,\vartheta-2\delta_2-\delta_3,y_2,w),0,
    \delta_{1:4}\big)+2\delta_5+\delta_6,\\
  f_2(y_1,y_2,\delta,w)
   &:=\mathcal{G}^*_2\big(\iota(y_1,0,y_2,w),0,\delta_{1:4}\big),\\
  g(w)&:=\frac{1}{k_dw_2},
  \end{aligned}
  \end{equation}
  where $\delta=(\delta_1,\dots,\delta_6)\in\mathbb{R}^6$ and
  $\delta_{1:4}:=(\delta_1,\dots,\delta_4)$.  These are well defined,
  since for $(y_1,y_2,\delta,w)\in\mathcal{V}$ and $s\in\mathbb{R}$
  we have $\iota(y_1,s,y_2,w)\in \mathcal{O}$.
\item[(iii)] \emph{Membership.}  The system so obtained belongs to
  $\mathcal{B}(\mathcal{Z})$, with
  $\gamma=1/(k_d\tau_E)\ge1/(k_d\tau_{\max})>0$ and with the growth
  constant $\bar B$ of \eqref{eq:Bbar}.
\end{enumerate}
\end{proposition}

\begin{proof}
\emph{(i), existence and uniqueness.}  Fix $\zeta=(y_1,y_2)\in\mathcal{Z}$
and abbreviate $q(t):=\mathcal{Q}(d_1(t),y_1(t),y_2(t))$, so that
$q(t)\in[\underline{\mathcal Q},\overline{\mathcal Q}]$ for almost all
$t\ge 0$.  The right-hand side of \eqref{eq:intdyn} is measurable in $t$,
since $d_1$ is measurable and $d_2,\zeta$ are continuous, and locally
Lipschitz in $\eta=(n_\alpha,\tau_E)\in\mathbb{R}\times(0,\infty)$ uniformly in $t$:
$\tau\mapsto q(t)\tau^{\alpha_P}-\tau$ is smooth on $(0,\infty)$ with
$q$ bounded, and $R_*$ is locally Lipschitz on $\mathbb{D}$ with
$(y_1(t),y_2(t))$ in a compact set.  Hence \eqref{eq:intdyn} has a
unique maximal solution $\eta:[0,\omega)\to\mathbb{R}\times(0,\infty)$, $0<\omega\le\infty$.

\emph{(i), the bounds \eqref{eq:bibobounds}.}  Since $\alpha_P\in(0,1)$
and $q(t)\in[\underline{\mathcal Q},\overline{\mathcal Q}]$ we have
\[
q(t)\tau^{\alpha_P}-\tau\ge0\ \text{ for }
 \tau\le\underline{\mathcal Q}^{\frac{1}{1-\alpha_P}},
\qquad
q(t)\tau^{\alpha_P}-\tau\le0\ \text{ for }
 \tau\ge\overline{\mathcal Q}^{\frac{1}{1-\alpha_P}},
\]
for almost all $t\ge 0$, so $\tau_E$, which is absolutely continuous, cannot leave
$[\tau_{\min},\tau_{\max}]$: this is the first inequality.  For the
second, $n_{DT}^*(\mathbf{y})\in[0,y_1]$ for every $\mathbf{y}=(y_1,y_2,a)\in\mathbb{D}$,
so \eqref{eq:Rext} and \eqref{eq:svmax} give
\[
0\le R_*\big(y_1(t),y_2(t),n_\alpha(t)\big)
 \le\tfrac14\overline y_1^2\sigma_v^{\max}
\qquad\text{for all }t\in[0,\omega),
\]
and therefore
\[
\tfrac12 \tfrac{\mathrm d}{\mathrm dt} n_\alpha(t)^2
 \le-\frac{n_\alpha(t)^2}{k_\alpha\tau_{\max}}+\Sigma |n_\alpha(t)|
\quad\text{a.e.},\qquad\text{hence}\qquad
|n_\alpha(t)|\le\max\{n_\alpha^0,\,k_\alpha\tau_{\max}\Sigma\}.
\]
Since $R_*\ge0$ and $d_2\ge0$ by (D2), moreover
$\dot n_\alpha(t)\ge-n_\alpha(t)/(k_\alpha\tau_{\min})$ a.e., so $n_\alpha(t)\ge n_\alpha^0 e^{-t/(k_\alpha\tau_{\min})}\ge0$. Both components are thus bounded on bounded intervals and $\tau_E$ is
bounded away from $0$, so $\omega=\infty$.

\emph{(i), property (B3).}  Causality follows from uniqueness of the solution on
$[0,t]$, which is thus determined by $\zeta|_{[0,t]}$, $d_1|_{[0,t]}$,
$d_2|_{[0,t]}$ and $\eta^0$.  Since the
right-hand side of \eqref{eq:intdyn} is locally Lipschitz in
$(\zeta,\eta)$ jointly, uniformly in $t$, local Lipschitz continuity of the operator $\cT$ follows from
Gronwall's inequality. It remains to show the bounded-input bounded-output property: the bounds
\eqref{eq:bibobounds} do not depend on $\zeta$, so with
$\mathcal{K}_\eta:=[-\max\{n_\alpha^0,k_\alpha\tau_{\max}\Sigma\},
\max\{n_\alpha^0,k_\alpha\tau_{\max}\Sigma\}]\times[\tau_{\min},\tau_{\max}]$
and $\mathcal{K}_d\subseteq\mathbb{R}^6$ a compact set containing the
essential image of the disturbance $d$ of (ii), which exists since
$d_1,\dots,d_4,\dot d_2,\dot d_3$ are essentially bounded by (D1)--(D3),
the set
\[
\mathcal{K}:=[\underline y_1,\overline y_1]\times
 [\underline y_2,\overline y_2]\times\mathcal{K}_d\times\mathcal{K}_\eta
\]
is a compact subset of $\mathcal{V}$ defined in~\eqref{eq:Vplant}, because
$\underline y_i>0$ and $\tau_{\min}>0$. Then we have $(\zeta(t),d(t),\mathcal{T}(\zeta)(t))\in\mathcal{K}$ for almost all $t\ge 0$, so (B3) is shown.

\emph{(ii), the diffeomorphism.}  $\Xi$ is smooth on
$\mathbb{R}^4\times(0,\infty)$ by definition of $\Lambda(x)$ in~\eqref{eq:Lambda}, and $\iota$ is
smooth there by \eqref{eq:nninv}.  For $x\in\mathbb{R}^4\times(0,\infty)$
one has $\iota(\Xi(x))=x$: the components $n_\alpha$, $E$, $\tau_E$ are
carried through, $n_{DT}=h_1(x)-2n_\alpha$, and $n_n$ is recovered by
solving $\mathcal{G}_1^0(x)=\Lambda(x)/\tau_E$ for $n_n$ in
\eqref{eq:Lambda}.  A similar computation gives $\Xi(\iota(y_1,z,y_2,n_\alpha,\tau_E))=(y_1,z,y_2,n_\alpha,\tau_E)$ for all
$(y_1,z,y_2,n_\alpha,\tau_E)\in \mathbb{R}^4\times(0,\infty)$.

\emph{(ii), the internal state.}  Let $x$ be a solution as in the
statement and put $\zeta:=(h_1(x),h_2(x))=(y_1,y_2)\in\mathcal{Z}$; all
differential identities below hold for almost all $t$.  Since
$x(t)\in\mathcal{X}$, \eqref{eq:Rext} gives
$R(t)=R_*(y_1(t),y_2(t),n_\alpha(t))$, so \eqref{eq:na} is the first
line of \eqref{eq:intdyn}; \eqref{eq:closure} with $n_e=y_1$ and $E=y_2$
is its second line; and $\eta(x(0))=\eta^0$ by hypothesis.  Hence
$\eta(x(\cdot))$ solves \eqref{eq:intdyn} with this $\zeta$, and by the
uniqueness in (i),
\begin{equation}\label{eq:etaT}
\eta\big(x(t)\big)=\mathcal{T}(\zeta)(t)\qquad\text{for all }t\ge0 .
\end{equation}

\emph{(ii), the class form.}  By \eqref{eq:Lambda},
\begin{equation}\label{eq:mu}
\dot y_1(t)=\mathcal{G}_1^0\big(x(t)\big)+\mu(t),\qquad
\mu(t):=2d_2(t)+d_3(t) ,
\end{equation}
so that, by \eqref{eq:etaT} and $\iota\circ\Xi=\mathrm{id}$,
\begin{equation}\label{eq:xiota}
x(t)=\iota\big(y_1(t),\,\dot y_1(t)-\mu(t),\,y_2(t),\,
 \mathcal{T}(\zeta)(t)\big).
\end{equation}
Differentiating \eqref{eq:mu}, which is legitimate because
$d_2,d_3\in W^{1,\infty}(\R_{\ge 0};\R)$ by (D2), gives \eqref{eq:ddoty1}; since
$x(t)\in\mathcal{X}$, $f$ may be replaced by $f^*$ there, and with the
affine structure \eqref{eq:affine} and the substitution
\eqref{eq:xiota} we obtain
\[
\begin{aligned}
\ddot y_1(t)&=\mathcal{G}^*_{11}\big(x(t),0,(d_1(t),\dots,d_4(t))\big)
 +2\dot d_2(t)+\dot d_3(t)+\frac{u_1(t)}{k_d\tau_E(t)}\\
 &=f_1\big(y_1(t),\dot y_1(t),y_2(t),d(t),\mathcal{T}(\zeta)(t)\big)
 +g\big(\mathcal{T}(\zeta)(t)\big)\,u_1(t),\\
\dot y_2(t)&=\mathcal{G}^*_2\big(x(t),0,(d_1(t),\dots,d_4(t))\big)+u_2(t)
 =f_2\big(y_1(t),y_2(t),d(t),\mathcal{T}(\zeta)(t)\big)+u_2(t),
\end{aligned}
\]
with $f_1,f_2,g$ from \eqref{eq:fdef} and $d$ the disturbance of (ii);
in $f_2$ the second argument of $\iota$ may be set to $0$ because
$\mathcal{G}^*_2=f^*_4$ does not involve $n_n$.  This is
\eqref{eq:class}, \eqref{eq:Fdef} with $a_2=1$.

\emph{(iii), properties (B1) and (B3).}  $g$ is continuous
and nonzero on $\mathcal{W}$, and $f_2$ is a
composition of $\iota$, the starred maps of \eqref{eq:Rext} and the
right-hand side of \eqref{eq:E}, all of which are continuous wherever
$y_1,y_2>0$ and $w_2>0$; this is (B1).
(B3) was shown in (i).  Along solutions, $\gamma(t)=1/(k_d\tau_E(t))\ge
1/(k_d\tau_{\max})>0$ by \eqref{eq:bibobounds}.

\emph{(iii), property (B2).}  Insert \eqref{eq:nninv} into $f_1$.  The
variable $\vartheta$ enters $f_1$ only through the second argument of
$\iota$, i.e.\ through $n_n$ and through
$\Lambda=w_2(\vartheta-2\delta_2-\delta_3)$, and
\begin{equation}\label{eq:df1}
\frac{\partial f_1}{\partial\vartheta}(y_1,\vartheta,y_2,\delta,w)
 =-\frac{1}{w_2}\Big(\frac{1}{k_d}+\frac{1}{k_{DT}}\Big)
  -\frac{\mathcal{Q}(\delta_1,y_1,y_2)\,w_2^{\alpha_P-1}-1}{T_c} ,
\end{equation}
the last term being $-\dot\tau_E/\tau_E$ evaluated at the point
$(y_1,y_2,\delta,w)$ by \eqref{eq:closure}; in particular $f_1$ is
affine in $\vartheta$, and the remaining dependence is continuous on
$\mathcal{V}_1$.  At every point
\[
 (y_1,y_2,\delta,w)\in   \mathcal{R} = \bigcup_{\zeta\in\mathcal{Z}}
   \operatorname{ess\,im}\big(\zeta,d,\mathcal{T}(\zeta)\big) \subseteq \cV
\]
one has
$w_2\in[\tau_{\min},\tau_{\max}]$ by \eqref{eq:bibobounds},
$\delta_1\in[\underline d_1,\overline d_1]$ since $\delta$ lies in the
essential image of $d$, and $y_i\in[\underline y_i,\overline y_i]$,
hence, as $\alpha_P-1<0$,
\[
\mathcal{Q}(\delta_1,y_1,y_2)\,w_2^{\alpha_P-1}\in
 \big[\underline{\mathcal Q}\,\tau_{\max}^{\alpha_P-1},\;
 \overline{\mathcal Q}\,\tau_{\min}^{\alpha_P-1}\big],
\]
and therefore, by \eqref{eq:df1} and \eqref{eq:Bbar},
\[
\Big|\frac{\partial f_1}{\partial\vartheta}\Big|
 \le\frac{1}{\tau_{\min}}\Big(\frac{1}{k_d}+\frac{1}{k_{DT}}\Big)
 +\frac{\Theta}{T_c}=\bar B .
\]  Hence (B2) holds with
$c(y_1,y_2,\delta,w):=|f_1(y_1,0,y_2,\delta,w)|$ and this $\bar B$.
\end{proof}

Membership is a property of the plant and the box alone.  It is
established for every continuous history with values in the box,
whether or not any controller produces it, and presupposes no
closed loop: the box is the domain on which the internal
dynamics is shown to be well posed and bounded, not a claim about where
the outputs will be.  The only condition the plant imposes on a box is
$\underline y_i>0$; whether a control design keeps the outputs inside
a given box is a question for that design, answered in the
next section by choosing the box from the design data and proving
invariance.

%=====================================================================

%=====================================================================
\section{Controller and main result}\label{sec:main}
%=====================================================================

\subsection{Input constraints}\label{sec:inputconstr}

Both actuators are constrained, $S(t)\in[0,S_{\max}]$ and
$P_{\mathrm{aux}}(t)\in[0,P_{\max}]$, with $S_{\max}$ the capacity of
the fuelling system and $P_{\max}$ the installed heating capacity.  From
now on the inputs of \eqref{eq:class} are the
relative deviations
\[
u_1:=S/\bar S-1,\qquad u_2:=P_{\mathrm{aux}}/\bar P-1
\]
from the nominal fuelling rate $\bar S$ and heating power density
$\bar P$ of the operating point of Section~\ref{sec:closure}.  This
affine change of input coordinates keeps the plant in the class: by
\eqref{eq:model}, $F_1+\gamma\bar S$, $\gamma\bar S$, $F_2+\bar P$ and
$\bar P$ take the places of $F_1,\gamma,F_2,a_2$ in
Proposition~\ref{prop:membership}, with the same operator, box and
growth constant $\bar B$, the shift of $F_1$ being absorbed by the
function $c$ of (B2).  The design needs the deviation form: the
demand of \eqref{eq:controller} below vanishes at zero error, so the
zero input must be the one that holds the operating point; and the
constraint sets must contain $0$ in their interior.  In these
variables the constraints are box constraints
\begin{equation}\label{eq:sat}
	\operatorname{sat}_i(v):=\min\{\max\{v,\underline u_i\},\overline u_i\},\quad i=1,2,
\end{equation}
for $\underline{u}_i<0<\overline{u}_i$, and we write
$M_i:=\max\{|\underline{u}_i|,\overline{u}_i\}$ and
$\underline{M}_i:=\min\{|\underline{u}_i|,\overline{u}_i\}>0$.
Of the saturation function \eqref{eq:sat} we use that it is Lipschitz, takes
values in $[\underline u_i,\overline u_i]$ and is the identity there, and
satisfies
\[
	\operatorname{sgn}\big(v-\operatorname{sat}_i(v)\big)=\operatorname{sgn}(v)
	\quad\text{whenever } v\ne\operatorname{sat}_i(v).
\]
The heating power
carries in addition a rate limit, $|\dot P_{\mathrm{aux}}(t)|\le\bar r$
with $\bar r$ of order $10$--$20\,\mathrm{MW\,s^{-1}}$ for ITER
\cite{SchKrsTyn03}, i.e.\ $|\dot u_2(t)|\le\hat r:=\bar r/(V\bar P)$,
enforced through
$\operatorname{sat}_{\hat r}(z):=\min\{\hat r,\max\{-\hat r,z\}\}$.

\subsection{The controller}

Consider the set of admissible funnels
\begin{equation*}
	\Phi:=\Big\{\varphi\in W^{1,\infty}(\mathbb{R}_{\ge0};\mathbb{R})
	\;\Big|\;\varphi>0,\ \inf_{t\ge0}\varphi(t)>0\Big\},
\end{equation*}
and, for $\varphi_i\in\Phi$, set $\psi_i(\cdot):=\varphi_i(\cdot)^{-1}$.  With
$y_{\mathrm{ref},1},y_{\mathrm{ref},2}\in W^{1,\infty}(\mathbb{R}_{\ge0};\mathbb{R})$ we propose
\begin{subequations}\label{eq:controller}
	\begin{align}
		e_1(t) &= y_1(t)-y_{\mathrm{ref},1}(t)+I_1(t), \label{eq:e1}\\
		k_1(t) &= \big(1-\varphi_1(t)^2e_1(t)^2\big)^{-1},\label{eq:k1}\\
		\theta(t) &= \xi(t)+k_1(t)e_1(t),\quad
		v_1(t) = -\theta(t)/\big(\hat\theta^2-\theta(t)^2\big),
		\label{eq:v1}\\
		\dot\xi(t) &= -\kappa_3\,\xi(t)+v_1(t), \label{eq:xi}\\
		\dot I_3(t) &= -\kappa_3\,I_3(t)
		+\nu_1\big(v_1(t)-\operatorname{sat}_1(v_1(t))\big),
		\label{eq:I2}\\
		\dot I_1(t) &= -\kappa_1\,I_1(t)+I_3(t), \label{eq:I1}\displaybreak[1]\\[3pt]
		e_2(t) &= y_2(t)-y_{\mathrm{ref},2}(t)+I_2(t), \label{eq:e2}\\
		k_2(t) &= \big(1-\varphi_2(t)^2e_2(t)^2\big)^{-1},\quad
		v_2(t) = -k_2(t)e_2(t), \label{eq:v2}\\
		e_3(t) &= u_2(t)-\operatorname{sat}_2\big(v_2(t)\big), \displaybreak[1]\label{eq:e3}\\
		\dot I_2(t) &= -\kappa_2\,I_2(t)
		+\nu_2\big(v_2(t)-\operatorname{sat}_2(v_2(t))\big),
		\label{eq:I0}\\
		\dot u_2(t) &= \operatorname{sat}_{\hat r}\big(-k_3\,e_3(t)\big)
		\label{eq:act}
	\end{align}
\end{subequations}
with design constants $\hat\theta>0$, $\kappa_1,\kappa_2,\kappa_3,\nu_1,\nu_2,k_3>0$
and initial
values $\xi(0)=\xi^0\in\R$, $I_j(0)=I_j^0\in\R$, $j=1,2,3$, and
$u_2(0)\in[\underline u_2,\overline u_2]$.  The plant receives
$u_1(t)=\operatorname{sat}_1(v_1(t))$ and the actuator state $u_2(t)$
of \eqref{eq:act}.  In the heating channel three signals are to be distinguished, the
demand $v_2$, the actuator target $\operatorname{sat}_2(v_2)$ and the applied input $u_2$:
\[
e_2\xrightarrow{\eqref{eq:v2}}\underbrace{v_2}_{\text{demand}}
\xrightarrow{\ \operatorname{sat}_2\ }
\underbrace{\operatorname{sat}_2(v_2)}_{\text{target}}
\xrightarrow{\eqref{eq:act}}
\underbrace{u_2}_{\text{applied}},
\qquad
v_2-\operatorname{sat}_2(v_2)\xrightarrow{\eqref{eq:I0}}I_2
\xrightarrow{\eqref{eq:e2}}e_2 .
\]
The correction $I_2$ responds to the deficit between demand and target, not to the lag $e_3$
between target and applied input, on which the rate limit acts; the rate limit can therefore be
active while $I_2\equiv0$.  Either way the applied input $u_i$ never leaves
$[\underline u_i,\overline u_i]$ for $i=1,2$ (for $u_2$ by Lemma~\ref{lem:act} below), so the sign constraints on $S$ and $P_{\mathrm{aux}}$ hold by
construction and the invariance of Remark~\ref{rem:physical} is in
force.

The signals $I_1,I_2,I_3$ are driven by the \emph{input deficit}
$v_i-\operatorname{sat}_i(v_i)$, which vanishes whenever the demanded
input is admissible, and enter the loop only by shifting the references
to $y_{\mathrm{ref},1}-I_1$ and $y_{\mathrm{ref},2}-I_2$; channel~1
needs two of them because it has relative degree two, $I_3$ being the
filtered and scaled deficit and $I_1$ the reference shift it drives.
The filter and $I_3$ share the rate $\kappa_3$, and $I_1$ and $I_2$
have their own rates $\kappa_1$ and $\kappa_2$.  The actuator \eqref{eq:act} enforces the rate
limit by construction and, by Lemma~\ref{lem:act}, the amplitude limit
with it; since $u_2$ thus stays in the constraint set, it needs no
correction term of its own.  The funnel boundaries $\psi_1,\psi_2$ and
$\hat\theta$ remain prescribed throughout; it is the reference, not the
funnel, that is adjusted due to saturation.

\subsection{Compatibility}\label{sec:compat}

Hard input constraints and prescribed output error bounds cannot both be met
for arbitrary reference signals and disturbances, so a compatibility condition is
unavoidable.  We ask that the correction states $I_2$ and $I_3$, which
are driven by the input deficits, stay bounded by constants fixed in
advance.  This is a condition on the closed-loop solution, see also
Remark~\ref{rem:cost}; in
Section~\ref{sec:apriori} it is reduced to conditions on the plant and
on the data, in the spirit of \cite[Thm.~1]{FotRov24}.

\emph{The design box.}  For numbers $\hat{I}_2\ge|I_2^0|$ and
$\hat{I}_3\ge|I_3^0|$ put
\begin{equation}\label{eq:I0hat}
	\hat{I}_1:=\max\{|I_1^0|,\hat{I}_3/\kappa_1\},
\end{equation}
for $i=1,2$, put
\begin{equation}\label{eq:Ypm}
	Y_i^+:=\|y_{\mathrm{ref},i}\|_\infty+\|\psi_i\|_\infty+\hat{I}_i,\qquad
	Y_i^-:=\inf_{t\ge0}y_{\mathrm{ref},i}(t)-\|\psi_i\|_\infty-\hat{I}_i ,
\end{equation}
and, for $\epsilon>0$,
\begin{equation}\label{eq:compatZ}
	\mathcal{Z}_\epsilon:=\big\{\zeta\in\mathcal{C}(\mathbb{R}_{\ge0};
	\mathbb{R}^2)\;\big|\;
	Y_i^--\epsilon\le \zeta_i(t)\le Y_i^++\epsilon\ \text{ for all }t\ge0,\
	i=1,2\big\}.
\end{equation}
All of these are formed from the design data, the initial values and
the two numbers $\hat{I}_2,\hat{I}_3$.  If the plant belongs to
$\mathcal{B}(\mathcal{Z}_\epsilon)$, the closed loop \eqref{eq:class},
\eqref{eq:controller} -- in which the plant receives
$u_1(t)=\operatorname{sat}_1(v_1(t))$ and the actuator state $u_2(t)$
of \eqref{eq:act} -- is defined, and a \emph{solution} of it on
$[0,\omega)$, $0<\omega\le\infty$, is a locally absolutely continuous
function $(y_1,\dot{y}_1,y_2,\xi,I_1,I_2,I_3,u_2):[0,\omega)\to
\mathbb{R}^8$ which satisfies \eqref{eq:class}, \eqref{eq:controller}
almost everywhere, the operator acting on the output history as
explained after Definition~\ref{def:class}, and which satisfies
$\varphi_i(t)|e_i(t)|<1$, $|\theta(t)|<\hat{\theta}$ and
$Y_i^--\epsilon<y_i(t)<Y_i^++\epsilon$ for all $t\in[0,\omega)$; it is
\emph{maximal} if it has no proper extension.  Since
$\mathcal{B}(\mathcal{Z})\subseteq\mathcal{B}(\mathcal{Z}')$ for
$\mathcal{Z}'\subseteq\mathcal{Z}$, membership in
$\mathcal{B}(\mathcal{Z})$ for any $\mathcal{Z}\supseteq
\mathcal{Z}_\epsilon$ suffices for membership in
$\mathcal{B}(\mathcal{Z}_\epsilon)$.

The following assumption is used in the main result to ensure existence of a solution for the closed-loop system such that the control objectives are achieved. It imposes a requirement on the closed loop itself, which is not desirable. Since $I_2$ and $I_3$ are states of the controller~\eqref{eq:controller}, the required bounds can be checked at every time instant along a solution, so that a violation is detected when it occurs. A sufficient condition, which excludes a violation in advance, is derived in Section~\ref{sec:apriori}.

\begin{assumption}\label{ass:compat}
	There are constants $\hat{I}_2\ge|I_2^0|$, $\hat{I}_3\ge|I_3^0|$ and
	$\epsilon>0$ such that the plant belongs to
	$\mathcal{B}(\mathcal{Z}_\epsilon)$ with the design box
	\eqref{eq:compatZ} formed from $\hat{I}_2,\hat{I}_3$, and such that every
	maximal solution of the closed loop \eqref{eq:class},
	\eqref{eq:controller} on its interval of existence $[0,\omega)$
	satisfies
	\[
	|I_2(t)|\le\hat{I}_2\quad\text{and}\quad|I_3(t)|\le\hat{I}_3
	\qquad\text{for all }t\in[0,\omega).
	\]
\end{assumption}

The condition is self-referential -- the constants enter the box on
which the closed loop is defined -- but not circular: for given
$\hat{I}_2,\hat{I}_3,\epsilon$ the box, the class, the closed loop and its
maximal solutions are defined, in this order, before any solution is
considered, and the two bounds are then a statement about those
solutions, namely that none of them, on however short an interval it
may exist, breaks them.

One thing is important to note at this point: That the plant~\eqref{eq:model}, \eqref{eq:closure}, \eqref{eq:io} belongs to the system class $\mathcal{B}(\mathcal{Z}_\epsilon)$ is ensured by  Proposition~\ref{prop:membership} exactly when
$Y_i^->\epsilon$, i.e.\ when funnel width and largest correction do not
reach down to a vanishing density or energy.  That the outputs of the
closed loop stay in the box is a conclusion of
Theorem~\ref{thm:main} below.

\subsection{Three preliminary lemmas}

\begin{lemma}\label{lem:cross}
	Let $J\subseteq\mathbb{R}$ be an interval containing its left endpoint
	$t_0:=\min J$, let $z:J\to\mathbb{R}$ be locally absolutely continuous
	and let $c\in\mathbb{R}$.  If $\dot{z}(t)\le0$ for almost all $t\in J$
	with $z(t)>c$, then $z(t)\le\max\{z(t_0),c\}$ for all $t\in J$.
\end{lemma}

\begin{proof}
	Put $c':=\max\{z(t_0),c\}$ and suppose $z(t_1)>c'$ for some
	$t_1\in J$.  Then $s:=\sup\{t\in[t_0,t_1]\mid z(t)\le c'\}$ satisfies
	$z(s)=c'$ by continuity, and $z>c'\ge c$ on $(s,t_1]$, so
	$\dot{z}\le0$ a.e.\ on $(s,t_1]$ and
	$z(t_1)=z(s)+\int_s^{t_1}\dot{z}(t)\,\mathrm dt\le c'$, a contradiction.
\end{proof}

\begin{lemma}[Actuator]\label{lem:act}
	Let $0<\omega\le\infty$, let $q:[0,\omega)\to[\underline{u}_2,\overline{u}_2]$
	be measurable, and let $u_2:[0,\omega)\to\mathbb{R}$ be locally
	absolutely continuous with $u_2(0)\in[\underline{u}_2,\overline{u}_2]$
	and
	\[
	\dot{u}_2(t)=\operatorname{sat}_{\hat{r}}\big(-k_3(u_2(t)-q(t))\big)
	\qquad\text{for almost all }t\in[0,\omega)
	\]
	(for every such $q$ this initial value problem has exactly one such
	solution $u_2$, since its right-hand side is measurable in $t$ and
	globally Lipschitz and bounded in $u_2$).  Then $e_3:=u_2-q$ satisfies
	\begin{equation}\label{eq:actbound}
		u_2(t)\in[\underline{u}_2,\overline{u}_2],\qquad
		|e_3(t)|\le W_2:=\overline{u}_2-\underline{u}_2
		\qquad\text{for all }t\in[0,\omega),
	\end{equation}
	and $|\dot{u}_2(t)|\le\hat{r}$ for almost all $t\in[0,\omega)$.
\end{lemma}

\begin{proof}
	For almost all $t$ with $u_2(t)>\overline{u}_2$ we have
	$u_2(t)>\overline{u}_2\ge q(t)$, hence $e_3(t)>0$ and, since
	$\operatorname{sat}_{\hat{r}}$ preserves signs,
	$\dot{u}_2(t)=\operatorname{sat}_{\hat{r}}(-k_3e_3(t))<0$.
	Lemma~\ref{lem:cross} with $z=u_2$ and $c=\overline{u}_2$ gives
	$u_2(t)\le\max\{u_2(0),\overline{u}_2\}=\overline{u}_2$ for all
	$t\in[0,\omega)$, and with $z=-u_2$ and $c=-\underline{u}_2$ it gives
	$u_2(t)\ge\underline{u}_2$ in the same way.  Since $u_2(t)$ and $q(t)$
	both lie in $[\underline{u}_2,\overline{u}_2]$, \eqref{eq:actbound}
	follows; the rate bound is immediate.
\end{proof}

\begin{lemma}[Shadow filter]\label{lem:shadow}
	Let $0<\omega\le\infty$, let $v_1:[0,\omega)\to\mathbb{R}$ be
	measurable and locally essentially bounded, and let
	$\xi,I_3:[0,\omega)\to\mathbb{R}$ be locally absolutely continuous
	with $\xi(0)=\xi^0$, $I_3(0)=I_3^0$ and \eqref{eq:xi}, \eqref{eq:I2}
	for almost all $t\in[0,\omega)$.  Let $\xi_s:[0,\omega)\to\mathbb{R}$
	be the solution of
	$\dot{\xi}_s(t)=-\kappa_3\,\xi_s(t)+\operatorname{sat}_1\big(v_1(t)\big)$,
	$\xi_s(0)=\xi^0-\nu_1^{-1}I_3^0$.  Then, for all $t\in[0,\omega)$,
	\begin{equation}\label{eq:shadow}
		\xi_s(t)=\xi(t)-\nu_1^{-1}I_3(t),\qquad
		|\xi_s(t)|\le\bar{\xi}_s
		:=\max\big\{|\xi^0-\nu_1^{-1}I_3^0|,\,M_1/\kappa_3\big\}.
	\end{equation}
\end{lemma}

\begin{proof}
	Put $\chi(t):=\xi(t)-\nu_1^{-1}I_3(t)$.  By \eqref{eq:xi} and
	\eqref{eq:I2}, for almost all $t$,
	\[
	\dot{\chi}(t)=-\kappa_3\,\xi(t)+v_1(t)+\nu_1^{-1}\kappa_3\,I_3(t)
	-\big(v_1(t)-\operatorname{sat}_1(v_1(t))\big)
	=-\kappa_3\,\chi(t)+\operatorname{sat}_1\big(v_1(t)\big),
	\]
	and $\chi(0)=\xi_s(0)$.  Thus $\chi$ and $\xi_s$ solve the same linear
	initial-value problem with the same locally integrable inhomogeneity,
	so they coincide.  Variation of constants and
	$|\operatorname{sat}_1(v_1(t))|\le M_1$ give
	\[
	|\xi_s(t)|\le\mathrm e^{-\kappa_3 t}|\xi_s(0)|
	+\big(1-\mathrm e^{-\kappa_3 t}\big)\frac{M_1}{\kappa_3} ,
	\]
	whose right-hand side is a convex combination of $|\xi_s(0)|$ and $M_1/\kappa_3$, hence \eqref{eq:shadow}.
\end{proof}

Along every solution of the closed loop, Lemma~\ref{lem:act} applies
with $q=\operatorname{sat}_2(v_2)$ and Lemma~\ref{lem:shadow} with the
solution's $v_1$, $\xi$ and $I_3$.

\subsection{Main result}

We are now in the position to state the main result.

\begin{theorem}\label{thm:main}
	Let $\varphi_1,\varphi_2\in\Phi$ and
	$y_{\mathrm{ref},1},y_{\mathrm{ref},2}\in W^{1,\infty}(\R_{\ge 0};\R)$, let the initial
	values in \eqref{eq:controller} satisfy $\varphi_1(0)|e_1(0)|<1$,
	$\varphi_2(0)|e_2(0)|<1$ and $|\theta(0)|<\hat{\theta}$, and let
	Assumption~\ref{ass:compat} hold with constants $\hat{I}_2,\hat{I}_3$,
	the box $\mathcal{Z}_\epsilon$ of \eqref{eq:compatZ} formed from them
	and $\hat{I}_1$ from \eqref{eq:I0hat}.  Then the closed-loop system
	\eqref{eq:class}, \eqref{eq:controller} has a solution, every solution
	extends to a maximal one, and every maximal solution is global and
	satisfies:
	\begin{enumerate}
		\item[(i)] $v_1,v_2,k_1,k_2,\xi,u_2,e_3$ and $\dot{y}_1$ are bounded,
		$|I_1(t)|\le\hat{I}_1$ for all $t\ge0$, and $|\dot{u}_2(t)|\le\hat{r}$
		for almost all $t\ge0$;
		\item[(ii)] there are $\varepsilon_1,\varepsilon_2>0$ and
		$\varepsilon_\theta\in(0,\hat{\theta}]$ with, for all $t\ge0$,
		\[
		|e_i(t)|\le\psi_i(t)-\varepsilon_i\ (i=1,2),\qquad
		|\theta(t)|\le\hat{\theta}-\varepsilon_\theta ;
		\]
		\item[(iii)] the tracking errors obey, for all $t\ge0$ and $i=1,2$,
		\[
		\big|y_i(t)-y_{\mathrm{ref},i}(t)\big|\le\psi_i(t)+|I_i(t)|
		\le\psi_i(t)+\hat{I}_i ;
		\]
		\item[(iv)] if $v_2(t)=\operatorname{sat}_2(v_2(t))$ for all
		$t\in[t_0,t_1)$, then
		$|I_2(t)|=\mathrm e^{-\kappa_2(t-t_0)}|I_2(t_0)|$ there; if
		$v_1(t)=\operatorname{sat}_1(v_1(t))$ for all $t\in[t_0,t_1)$, then there
		\begin{align*}
		|I_3(t)|&=\mathrm e^{-\kappa_3(t-t_0)}|I_3(t_0)|,\\
		|I_1(t)|&\le\big(|I_1(t_0)|+(t-t_0)|I_3(t_0)|\big)\,
		\mathrm e^{-\min\{\kappa_1,\kappa_3\} (t-t_0)}.
		\end{align*}
	\end{enumerate}
\end{theorem}

\begin{proof}
	\emph{Step 1 (existence).}  Consider the relatively open set
	\begin{multline*}
		\mathcal{D}:=\bigg\{(t,y_1,\hat{y}_1,y_2,\xi,I_1,I_2,I_3,u_2)
		\in\mathbb{R}_{\ge0}\times\mathbb{R}^8\ \big|\
		\varphi_i(t)|y_i-y_{\mathrm{ref},i}(t)+I_i|<1,\\
		\bigg|\xi+\frac{y_1-y_{\mathrm{ref},1}(t)+I_1}{1-\varphi_1(t)^2 (y_1-y_{\mathrm{ref},1}(t)+I_1)^2}\bigg|<\hat{\theta},\ Y_i^--\epsilon<y_i<Y_i^++\epsilon\ (i=1,2)
		\bigg\},
	\end{multline*}
	where $Y_i^\pm$, $\epsilon$ are those of
	Assumption~\ref{ass:compat}.  The solutions of~\eqref{eq:class},~\eqref{eq:controller} in the
	sense of Section~\ref{sec:compat} are exactly the solutions whose graph evolves
	in $\mathcal{D}$, and the initial point lies in $\mathcal{D}$ by the
	hypotheses on $e_1(0)$, $e_2(0)$, $\theta(0)$ and since
	$|I_1^0|\le\hat{I}_1$, $|I_2^0|\le\hat{I}_2$.  On $\mathcal{D}$ the
	right-hand side of the differential equations in~\eqref{eq:controller} is measurable in $t$ and
	locally Lipschitz in the remaining arguments, and every output history
	with graph in $\mathcal{D}$ lies in $\mathcal{Z}_\epsilon$, where the
	operator $\cT$ from~\eqref{eq:class} is defined.  Since \cite[Thm.~B.1]{IlcRya09} is stated for
	operators on all of $\mathcal{C}(\mathbb{R}_{\ge0};\mathbb{R}^2)$, let
	$\Pi:\mathbb{R}^2\to\mathbb{R}^2$ be the projection onto the box
	$[Y_1^--\epsilon,Y_1^++\epsilon]\times[Y_2^--\epsilon,Y_2^++\epsilon]$;
	it is Lipschitz and the identity on the box, so
	$\zeta\mapsto\mathcal{T}(\Pi\circ\zeta)$ is an operator on
	$\mathcal{C}(\mathbb{R}_{\ge0};\mathbb{R}^2)$ which is causal, locally
	Lipschitz and bounded-input bounded-output and coincides with
	$\mathcal{T}$ on $\mathcal{Z}_\epsilon$.  By (B1) and (B2), $f_1,f_2,g$
	are continuous, so \cite[Thm.~B.1]{IlcRya09} applies: a solution
	exists, every solution extends to a maximal one, and for every maximal
	solution $(y_1,\dot y_1,y_2,\xi,I_1,I_2,I_3,u_2):[0,\omega)\to\R^8$, the closure of the
	graph is not a compact subset of $\mathcal{D}$.

	\emph{Step 2 (a priori bounds).}  On $[0,\omega)$ we have
	$|e_i(t)|<\psi_i(t)\le\|\psi_i\|_\infty$, $i=1,2$, and $|\theta(t)|<\hat{\theta}$ by the
	definition of $\mathcal{D}$, $|u_1(t)|\le M_1$ by \eqref{eq:sat},
	$|u_2(t)|\le M_2$ by Lemma~\ref{lem:act}, and $|I_2(t)|\le\hat{I}_2$,
	$|I_3(t)|\le\hat{I}_3$ by Assumption~\ref{ass:compat}. Let $\xi_s:[0,\omega)\to\mathbb{R}$
	be the solution of
	$\dot{\xi}_s(t)=-\kappa_3\,\xi_s(t)+\operatorname{sat}_1\big(v_1(t)\big)$,
	$\xi_s(0)=\xi^0-\nu_1^{-1}I_3^0$, then	$|\xi_s(t)|\le\bar{\xi}_s$ by
	Lemma~\ref{lem:shadow} for $\bar \xi_s$ defined in~\eqref{eq:shadow}. Variation of
	constants in \eqref{eq:I1} gives
	\[
	|I_1(t)|\le\mathrm e^{-\kappa_1t}|I_1^0|
	+\big(1-\mathrm e^{-\kappa_1t}\big)\frac{\hat{I}_3}{\kappa_1}
	\stackrel{\eqref{eq:I0hat}}{\le}\hat{I}_1 ,
	\]
	so that, by \eqref{eq:e1} and \eqref{eq:e2}, and invoking the definition of $Y_i^\pm$ in~\eqref{eq:Ypm},
	\[
	Y_i^-\le y_i(t)\le Y_i^+\qquad\text{for all }t\in[0,\omega),\ i=1,2 .
	\]
	By (B3) the arguments of $f_2$ and $g$ stay in a compact subset
	$\mathcal{K}\subseteq\mathcal{V}$; by (B1) and (B2), $F_2$, $\gamma$
	and the function $c$ of the growth bound are therefore bounded on
	$[0,\omega)$, say by $M_2^F$, $\Gamma_{\max}$ and $M_1^F$, so that
	\[
	|F_1(t)|\le M_1^F+\bar{B}|\dot{y}_1(t)|,\quad
	|F_2(t)|\le M_2^F,\quad|\gamma(t)|\le\Gamma_{\max}
	\quad\text{for almost all $t\in[0,\omega)$.}
	\]

	\emph{Step 3 (boundedness of $k_1$).}  By \eqref{eq:v1},
	$k_1(t)e_1(t)=\theta(t)-\xi(t)$, and by Lemma~\ref{lem:shadow} and
	Step 2, $|\xi(t)|\le\bar{\xi}_s+\nu_1^{-1}\hat{I}_3$.  With
	$|\theta(t)|<\hat{\theta}$,
	\begin{equation}\label{eq:A1}
		k_1(t)|e_1(t)|\le\hat{\theta}+\bar{\xi}_s+\nu_1^{-1}\hat{I}_3
		\qquad\text{for all }t\in [0,\omega).
	\end{equation}
	Since $k_1=1+\varphi_1^2e_1^2k_1$ by \eqref{eq:k1} and
	$\varphi_1|e_1|<1$, \eqref{eq:A1} gives, for all $t\in[0,\omega)$,
	\begin{equation}\label{eq:e1bound}
		k_1(t)\le1+\varphi_1(t)k_1(t)|e_1(t)|
		\le1+\|\varphi_1\|_\infty\big(\hat{\theta}+\bar{\xi}_s+\nu_1^{-1}\hat{I}_3\big) ,
	\end{equation}
	so $k_1\in L^\infty([0,\omega);\mathbb{R})$.  By \eqref{eq:k1},
	$\varphi_1(t)^2e_1(t)^2=1-k_1(t)^{-1}\le1-\|k_1\|_\infty^{-1}$, hence
	$|e_1(t)|\le\varpi_1^\ast\psi_1(t)\le\psi_1(t)-\varepsilon_1$ for all
	$t\in[0,\omega)$, where
	$\varpi_1^\ast:=\big(1-\|k_1\|_\infty^{-1}\big)^{1/2}<1$ and
	$\varepsilon_1:=(1-\varpi_1^\ast) \inf_{\tau\ge0}\psi_1(\tau)>0$.

	\emph{Step 4 (boundedness of $\dot{y}_1$).}  By
	\eqref{eq:class}, Step 2 and $|u_1|\le M_1$,
	\begin{equation}\label{eq:pdot}
		|\ddot y_1(t)|\le M_1^F+\Gamma_{\max}M_1+\bar{B}|\dot y_1(t)|
		\qquad\text{for almost all }t\in[0,\omega).
	\end{equation}
	Put $\lambda:=\max\{|\dot{y}_1(0)|,2Y_1^+\}+M_1^F+\Gamma_{\max}M_1$ and
	fix $t\in[0,\omega)$.  There is $s\in[0,t]$ with $t-s\le1$ and
	$|\dot{y}_1(s)|\le\max\{|\dot{y}_1(0)|,2Y_1^+\}$:
	if $t\le1$, take $s:=0$; if $t>1$, then, since $y_1$ is continuously
	differentiable, the mean value theorem provides $s\in(t-1,t)$ with
	$\dot{y}_1(s)=y_1(t)-y_1(t-1)$, and $|y_1|\le Y_1^+$ by Step 2 (note
	$-Y_1^+\le Y_1^-$).  Since $\dot{y}_1$ is locally absolutely
	continuous and $r-s\le1$ for $r\in[s,t]$, \eqref{eq:pdot}
	gives, for all $r\in[s,t]$,
	\[
	|\dot{y}_1(r)|\le|\dot{y}_1(s)|+\int_s^r|\ddot{y}_1(\tau)|\,\mathrm d\tau
	\le\lambda+\bar{B}\int_s^r|\dot{y}_1(\tau)|\,\mathrm d\tau .
	\]
	Then Gronwall's inequality gives
	$|\dot{y}_1(r)|\le\lambda\,\mathrm e^{\bar{B}(r-s)}$ for all
	$r\in[s,t]$.  In particular,
	$|\dot{y}_1(t)|\le\lambda\,\mathrm e^{\bar{B}(t-s)}
	\le\lambda\,\mathrm e^{\bar{B}}$, and since $t\in[0,\omega)$ was
	arbitrary,
	\begin{equation}\label{eq:pbound}
		\|\dot{y}_1\|_\infty\le\lambda\,\mathrm e^{\bar{B}} .
	\end{equation}
	By \eqref{eq:e1}, \eqref{eq:I1} and Step 2,
	$\dot{e}_1=\dot{y}_1-\dot{y}_{\mathrm{ref},1}-\kappa_1I_1+I_3$ is
	essentially bounded as well.

	\emph{Step 5 (boundedness of $k_2$).}  By \eqref{eq:class}, \eqref{eq:e2},
	\eqref{eq:I0} and $v_2=-k_2e_2$, for almost all $t\in[0,\omega)$,
	\[
	\begin{aligned}
		\dot{e}_2(t)&=-\nu_2k_2(t)e_2(t)+\tilde{g}_2(t),\\
		\tilde{g}_2(t)&:=F_2(t)-\dot{y}_{\mathrm{ref},2}(t)+a_2u_2(t)
		-\kappa_2I_2(t)-\nu_2\operatorname{sat}_2\big(v_2(t)\big),
	\end{aligned}
	\]
	and $\tilde{g}_2$ is essentially bounded by Step 2.
	Put $\varpi_2(t):=\varphi_2(t)|e_2(t)|\in[0,1)$.  The function
	$\varpi_2^2=\varphi_2^2e_2^2$ is locally absolutely continuous, and
	since $\varphi_2^2k_2e_2^2=\varpi_2^2/(1-\varpi_2^2)$ and $\varpi_2\le1$,
	\[
	\tfrac12\tfrac{\mathrm d}{\mathrm dt}\varpi_2(t)^2
	=\frac{\dot{\varphi}_2(t)}{\varphi_2(t)}\varpi_2(t)^2
	+\varphi_2(t)^2e_2(t)\dot{e}_2(t)
	\le\varpi_2(t)\Big(C_2-\nu_2\frac{\varpi_2(t)}{1-\varpi_2(t)^2}\Big)
	\]
	for almost all $t\in[0,\omega)$, where
	$C_2:=\|\dot{\varphi}_2/\varphi_2\|_\infty
	+\|\varphi_2\|_\infty\|\tilde{g}_2\|_\infty$.  As
	$\varpi\mapsto\nu_2\varpi/(1-\varpi^2)$ is an increasing bijection of
	$[0,1)$ onto $[0,\infty)$, there is a unique
	$\bar{\varpi}_2\in[0,1)$ with
	$\nu_2\bar{\varpi}_2/(1-\bar{\varpi}_2^2)=C_2$, and
	$\tfrac{\mathrm d}{\mathrm dt}\varpi_2^2\le0$ almost everywhere on
	$\{t\in[0,\omega)\mid\varpi_2(t)>\bar{\varpi}_2\}$.  Lemma~\ref{lem:cross} with
	$z=\varpi_2^2$ and $c=\bar{\varpi}_2^2$ gives
	$\varpi_2(t)\le\varpi_2^\ast:=\max\{\varpi_2(0),\bar{\varpi}_2\}<1$ for all $t\in [0,\omega)$,
	i.e.\ $|e_2(t)|\le\psi_2(t)-\varepsilon_2$ with
	$\varepsilon_2:=(1-\varpi_2^\ast) \inf_{\tau\ge0}\psi_2(\tau)>0$, and
	$k_2=(1-\varpi_2^2)^{-1}$ is bounded.

	\emph{Step 6 (boundedness of $v_1$).}  By \eqref{eq:xi} and
	$\xi=\theta-k_1e_1$, for almost all $t\in[0,\omega)$,
	\[
		\dot{\theta}(t)=-\kappa_3\,\theta(t)+B(t)+v_1(t),\qquad
		B(t):=\kappa_3\,k_1(t)e_1(t)+\tfrac{\mathrm d}{\mathrm dt}
		\big(k_1(t)e_1(t)\big).
	\]
	Here $\dot{k}_1=2k_1^2\big(\varphi_1\dot{\varphi}_1e_1^2
	+\varphi_1^2e_1\dot{e}_1\big)$, so that, with $\varphi_1|e_1|<1$ and
	$|e_1|\le\|\psi_1\|_\infty$, for almost all $t\in[0,\omega)$,
	\begin{equation}\label{eq:MB}
		|B(t)|\le\kappa_3\|k_1\|_\infty\|\psi_1\|_\infty
		+\|k_1\|_\infty\|\dot{e}_1\|_\infty
		+2\|k_1\|_\infty^2\big(\|\dot{\varphi}_1\|_\infty\|\psi_1\|_\infty^2
		+\|\dot{e}_1\|_\infty\big).
	\end{equation}
	Put $\varpi_\theta(t):=|\theta(t)|/\hat{\theta}\in[0,1)$.  The function
	$\varpi_\theta^2=\theta^2/\hat{\theta}^2$ is locally absolutely continuous, and
	since $\theta v_1=-\varpi_\theta^2/(1-\varpi_\theta^2)$ by \eqref{eq:v1},
	\[
	\tfrac12\tfrac{\mathrm d}{\mathrm dt}\varpi_\theta(t)^2
	=\frac{-\kappa_3\theta(t)^2+\theta(t)B(t)+\theta(t)v_1(t)}{\hat{\theta}^2}
	\le\frac{\varpi_\theta(t)}{\hat{\theta}^2}
	\Big(C_\theta-\frac{\varpi_\theta(t)}{1-\varpi_\theta(t)^2}\Big)
	\]
	for almost all $t\in[0,\omega)$, where
	$C_\theta:=\hat{\theta}\|B\|_\infty$.  As
	$\varpi\mapsto\varpi/(1-\varpi^2)$ is an increasing bijection of
	$[0,1)$ onto $[0,\infty)$, there is a unique
	$\bar{\varpi}_\theta\in[0,1)$ with
	$\bar{\varpi}_\theta/(1-\bar{\varpi}_\theta^2)=C_\theta$, and
	$\tfrac{\mathrm d}{\mathrm dt}\varpi_\theta^2\le0$ almost everywhere on
	$\{t\in[0,\omega)\mid\varpi_\theta(t)>\bar{\varpi}_\theta\}$.  Lemma~\ref{lem:cross} with
	$z=\varpi_\theta^2$ and $c=\bar{\varpi}_\theta^2$ gives
	$\varpi_\theta(t)\le\varpi_\theta^\ast:=\max\{\varpi_\theta(0),\bar{\varpi}_\theta\}<1$ for all $t\in [0,\omega)$,
	i.e.\ $|\theta(t)|\le\hat{\theta}-\varepsilon_\theta$ with
	$\varepsilon_\theta:=(1-\varpi_\theta^\ast)\hat{\theta}>0$, and, since
	$|v_1|=\hat{\theta}^{-1}\varpi_\theta/(1-\varpi_\theta^2)$ by \eqref{eq:v1},
	\begin{equation}\label{eq:v1bar}
		|v_1(t)|\le\max\big\{|v_1(0)|,\|B\|_\infty\big\}
		\qquad\text{for all }t\in[0,\omega).
	\end{equation}

	\emph{Step 7 ($\omega=\infty$ and assertions (i)--(iv)).}  Suppose
	$\omega<\infty$. By Steps 3 and 5, $\xi=\theta-k_1e_1$ and $v_2=-k_2e_2$ are bounded on $[0,\omega)$.
	By Step~2, \eqref{eq:pbound} and
	Lemma~\ref{lem:act}, the components $y_1,y_2$, $\dot{y}_1$, $I_1$ and $u_2$ of the solution are bounded on $[0,\omega)$, and $I_2$, $I_3$ are bounded by Assumption~\ref{ass:compat}.
	Furthermore, we have, for all $t\in [0,\omega)$,
	\[
	\varphi_i(t)|e_i(t)|\le\varpi_i^\ast<1,\qquad
	|\theta(t)|\le\varpi_\theta^\ast\hat{\theta}<\hat{\theta},\qquad
	Y_i^-\le y_i(t)\le Y_i^+,\qquad i=1,2,
	\]
	by Steps 2, 3, 5 and 6.  The closure of the graph of the
	solution is therefore a compact subset of $\mathcal{D}$, contradicting
	Step 1.  Hence $\omega=\infty$.  Assertion (i) then follows from Step~2, Lemma~\ref{lem:act} and the boundedness of the solution components shown above.
    Assertion~(ii) follows from Steps 3, 5 and 6, and Assertion~(iii) from \eqref{eq:e1}, \eqref{eq:e2} and $|e_i|<\psi_i$.  For Assertion~(iv), if $v_2=\operatorname{sat}_2(v_2)$ on
	$[t_0,t_1)$, then $\dot{I}_2=-\kappa_2I_2$ there by \eqref{eq:I0}, so $|I_2(t)| = \mathrm e^{-\kappa_2(t-t_0)}|I_2(t_0)|$.  If
	$v_1=\operatorname{sat}_1(v_1)$ on $[t_0,t_1)$, then
	$\dot{I}_3=-\kappa_3I_3$ and $\dot{I}_1=-\kappa_1I_1+I_3$ there by
	\eqref{eq:I2}, \eqref{eq:I1}, so $|I_3(t)| = \mathrm e^{-\kappa_3(t-t_0)}|I_3(t_0)|$ and
	\[
	\begin{aligned}
		|I_1(t)|&=\Big|\mathrm e^{-\kappa_1 (t-t_0)}I_1(t_0)
		+I_3(t_0) \int_0^{t-t_0} \mathrm e^{-\kappa_1(t-t_0-r)}\mathrm e^{-\kappa_3r}\,\mathrm dr\,
		\Big|\\
		&\le\big(|I_1(t_0)|+(t-t_0)|I_3(t_0)|\big)\,
		\mathrm e^{-\min\{\kappa_1,\kappa_3\}(t-t_0)}.
	\end{aligned}
	\]
\end{proof}

\begin{remark}\label{rem:cost}	
	Compared with the funnel control design involving a filter presented in~\cite{DenSchBer25},
	the input gain is not constant, but state-dependent here. Further note that funnel control is usually stated with a sign condition on the input
	gain, and Definition~\ref{def:class} has none. The reason is the
	input constraint: $|u_1|\le M_1$ bounds $\dot{y}_1$ through
	\eqref{eq:pbound} whatever the gain does, and the barrier on $\theta$
	is enforced through the filter \eqref{eq:xi}, into which $v_1$ enters
	with gain one.  Invariance of the funnels is thus a statement about
	the controller and the constraints, not about the plant.  The plant
	re-enters through Assumption~\ref{ass:compat}: compatibility can fail when the applied input drives the
	output in an adverse direction.  The invariance is thus conditional and, being independent of the sign of the gain,
	establishes tracking authority for neither sign.
	
	Finally, evolution of the tracking errors $y_i-y_{{\rm ref}, i}$ within the performance funnels defined by $\psi_i$
	is not achieved~-- no prescribed margin can be guaranteed unconditionally under hard input constraints and unknown disturbances~-- but the evolution of the \textit{corrected}
	tracking error $y_i-y_{{\rm ref}, i}+I_i$ within these bounds can be prescribed. An alternative, widening the funnel itself while
	the input saturates, is pursued in \cite{Ber24,BerImp26}.	
\end{remark}

%=====================================================================
\section{A priori verifiable conditions}\label{sec:apriori}
%=====================================================================

Assumption~\ref{ass:compat} is a hypothesis on the closed-loop
solution.  Following \cite[Thm.~1]{FotRov24}, here we replace it
by two conditions, one for each of the bounds on $I_2$ and $I_3$, in
terms of the data and of bounds on the plant over the design box, which
the controller does not use.

Write $M_2^{+}:=\overline{u}_2$, $M_2^{-}:=|\underline{u}_2|$, and, for $t\ge0$, $q^{\pm}(t):=\big(\psi_2(t)^2+4(M_2^{\pm})^2\big)^{1/2}$ and $\rho_2^{\pm}(t):=\frac{2M_2^{\pm}\psi_2(t)}{\psi_2(t)+q^{\pm}(t)}$
for the unique solution of $\rho/(1-\varphi_2(t)^2\rho^2)=M_2^{\pm}$ in
$[0,\psi_2(t))$; since $v_2(t)=-k_2(t)e_2(t)$, we have
$v_2(t)>\overline{u}_2$ exactly when $-e_2(t)>\rho_2^+(t)$, and
$v_2(t)<\underline{u}_2$ exactly when $e_2(t)>\rho_2^-(t)$.  Moreover,
$\dot{\rho}_2^{\pm}(t)=c^{\pm}(t)\dot{\psi}_2(t)$ almost everywhere with
$c^{\pm}(t):=8(M_2^{\pm})^3/\big(q^{\pm}(t)(\psi_2(t)+q^{\pm}(t))^2\big)\in(0,1]$,
so that
\begin{equation}\label{eq:rhoprop}
	0<\rho_2^{\pm}(t)<\psi_2(t)\qquad\text{and}\qquad
	\dot{\rho}_2^{\pm}(t)\le\|\dot{\psi}_2\|_\infty
	\ \text{ for almost all $t\ge 0$.}
\end{equation}
Finally, with $W_2$ from \eqref{eq:actbound}, put
$\bar{W}_2:=W_2+\hat{r}/k_3$.

\begin{proposition}\label{prop:ap}
	Let $\varphi_i$, $y_{\mathrm{ref},i}$, $i=1,2$, and the initial
	values be as in Theorem~\ref{thm:main}.  Let
	\begin{equation}\label{eq:ap2bound}
		\hat{I}_2:=a_2\bar{W}_2^{2}/(2\hat{r})+\|\psi_2\|_\infty
		+\max\Big\{|I_2^0|,\ (a_2\bar{W}_2
		+\|\dot{\psi}_2\|_\infty)/\kappa_2\Big\} ,
	\end{equation}
	let $\hat{I}_3\ge|I_3^0|$ and $\epsilon>0$, and let the plant belong to
	$\mathcal{B}(\mathcal{Z}_\epsilon)$ for the design box
	\eqref{eq:compatZ} formed from $\hat{I}_2,\hat{I}_3$.  Let
	$M_1^F,M_2^F,\Gamma_{\max}$ denote the suprema of $c$, $|F_2|$ and
	$|\gamma|$ over all output histories in $\mathcal{Z}_\epsilon$ and all
	$t\ge0$, with $c$ the function of (B2), and put
	\begin{equation}\label{eq:chain}
		\begin{aligned}
			\bar{k}_1&:=1+\|\varphi_1\|_\infty\big(\hat{\theta}+\bar{\xi}_s
			+\nu_1^{-1}\hat{I}_3\big),\\
			\hat{P}_1&:=\mathrm e^{\bar{B}}\big(\max\{|\dot{y}_1(0)|,2Y_1^+\}
			+M_1^F+\Gamma_{\max}M_1\big)
			+\|\dot{y}_{\mathrm{ref},1}\|_\infty+\kappa_1\hat{I}_1+\hat{I}_3,\\
			\bar{v}_1&:=\max\big\{|v_1(0)|,\ \kappa_3\bar{k}_1\|\psi_1\|_\infty
			+\bar{k}_1\hat{P}_1+2\bar{k}_1^2\big(\|\dot{\varphi}_1\|_\infty
			\|\psi_1\|_\infty^2+\hat{P}_1\big)\big\},
		\end{aligned}
	\end{equation}
	with $\bar{\xi}_s$ from \eqref{eq:shadow}, $Y_1^+$ from \eqref{eq:Ypm},
	$\hat{I}_1$ from \eqref{eq:I0hat}, and $v_1(0)$ given by
	\eqref{eq:e1}--\eqref{eq:v1} in terms of the initial values.  If
	\begin{align}
		M_2^F+\|\dot{y}_{\mathrm{ref},2}\|_\infty
		&\le a_2\underline{M}_2 ,\label{eq:ap2}\\
		\bar{v}_1&<\underline{M}_1+\kappa_3\hat{I}_3/\nu_1 ,\label{eq:ap1}
	\end{align}
	then Assumption~\ref{ass:compat} holds with $\hat{I}_2,\hat{I}_3$, and
	the conclusions of Theorem~\ref{thm:main} follow.
\end{proposition}

\begin{proof}
	Consider a solution $(y_1,\dot y_1,y_2,\xi,I_1,I_2,I_3,u_2):[0,\omega)\to\R^8$ of the closed-loop system
	\eqref{eq:class}, \eqref{eq:controller}.  By the
	definition of a solution in Section~\ref{sec:compat} its output
	history stays in the box of $\mathcal{Z}_\epsilon$, so,
	since $\cT$ is causal, for almost all $t\in[0,\omega)$,
	\begin{equation}\label{eq:plantb}
		|F_1(t)|\le M_1^F+\bar{B}|\dot{y}_1(t)|,\qquad
		|F_2(t)|\le M_2^F,\qquad|\gamma(t)|\le\Gamma_{\max} .
	\end{equation}

	\emph{Step 1 ($|I_2|\le\hat{I}_2$ on $[0,\omega)$).}  We show
	$I_2\le\hat{I}_2$.  Let $\ell:=\overline{u}_2-u_2$ and
	$h(\lambda):=(\lambda+\hat{r}/k_3)^2/(2\hat{r})$ for $\lambda\ge0$,
	and define
	\[
		V(t):=I_2(t)
		+\max\big\{{-}e_2(t),\rho_2^{+}(t)\big\}
		+a_2h\big(\ell(t)\big),\qquad t\in[0,\omega).
	\]
	By Lemma~\ref{lem:act}, $\ell$ takes values in $[0,W_2]$ and
	$|\dot{\ell}|\le\hat{r}$ almost everywhere.  Hence $V$ is locally
	absolutely continuous and $I_2\le V$; moreover
	$V(0)\le|I_2^0|+\psi_2(0)+a_2h(W_2)\le\hat{I}_2$ by
	\eqref{eq:rhoprop}, $|e_2(0)|<\psi_2(0)$ and
	$h(W_2)=\bar{W}_2^2/(2\hat{r})$.  By the choice of $\rho_2^+$, the
	relatively open set
	$\Sigma:=\{t\in[0,\omega)\mid v_2(t)>\overline{u}_2\}$ equals
	$\{t\in[0,\omega)\mid -e_2(t)>\rho_2^+(t)\}$.

	Off $\Sigma$ we have $v_2-\operatorname{sat}_2(v_2)\le0$, hence
	$\dot{I}_2\le-\kappa_2I_2$ by \eqref{eq:I0}, and
	$\max\{-e_2,\rho_2^+\}-\rho_2^+$ attains its minimum $0$, hence its
	derivative vanishes wherever it exists.  Moreover
	$\tfrac{\mathrm d}{{\mathrm d}t}h(\ell)=h'(\ell)\dot{\ell}
	\le\hat{r}h'(\ell)=\ell+\hat{r}/k_3\le\bar{W}_2$ and
	$I_2=V-\rho_2^+-a_2h(\ell)\ge V-\|\psi_2\|_\infty-a_2h(W_2)$.  With
	\eqref{eq:rhoprop} and \eqref{eq:ap2bound} this gives, almost
	everywhere on $[0,\omega)\setminus\Sigma$,
	\begin{equation}\label{eq:Voff}
		\dot{V}\le-\kappa_2\big(V-\|\psi_2\|_\infty-a_2h(W_2)\big)
		+\|\dot{\psi}_2\|_\infty+a_2\bar{W}_2
		\le-\kappa_2\big(V-\hat{I}_2\big) .
	\end{equation}

	On $\Sigma$ we have $\operatorname{sat}_2(v_2)=\overline{u}_2$, hence
	$e_3=-\ell$, so that \eqref{eq:act} reads
	$\dot{\ell}=-\min\{\hat{r},k_3\ell\}$, and
	$\tfrac{\mathrm d}{{\mathrm d}t}h(\ell)
	=-(\ell/\hat{r}+1/k_3)\min\{\hat{r},k_3\ell\}\le-\ell$
	because $\min\{a,b\}\ge ab/(a+b)$ for $a>0$, $b\ge0$.  Moreover,
	$V=-(y_2-y_{\mathrm{ref},2})+a_2h(\ell)$ there, by \eqref{eq:e2}.
	With \eqref{eq:class}, $u_2+\ell=\overline{u}_2\ge\underline{M}_2$,
	\eqref{eq:plantb} and \eqref{eq:ap2} we obtain, almost everywhere on
	$\Sigma$,
	\begin{equation}\label{eq:Von}
		\dot{V}\le-\big(F_2-\dot{y}_{\mathrm{ref},2}\big)-a_2u_2-a_2\ell
		\le M_2^F+\|\dot{y}_{\mathrm{ref},2}\|_\infty-a_2\underline{M}_2
		\le0 .
	\end{equation}

	By \eqref{eq:Voff} and \eqref{eq:Von}, $\dot{V}(t)\le0$ for
	almost all $t\in[0,\omega)$ with $V(t)>\hat{I}_2$, so
	Lemma~\ref{lem:cross} gives
	$I_2\le V\le\max\{V(0),\hat{I}_2\}=\hat{I}_2$ on $[0,\omega)$.  The
	bound $-I_2\le\hat{I}_2$ follows in the same way, with $-I_2$, $e_2$,
	$\rho_2^-$, $u_2-\underline{u}_2$ and
	$\{t\in[0,\omega)\mid v_2(t)<\underline{u}_2\}$ in place of
	$I_2$, $-e_2$, $\rho_2^+$, $\ell$ and $\Sigma$.

	\emph{Step 2 ($|I_3|\le\hat{I}_3$ on $[0,\omega)$).}  Let
	$w_1:=v_1-\operatorname{sat}_1(v_1)$ and
	\[
	\mathcal J:=\big\{t\in[0,\omega)\;\big|\;
	|I_3(r)|\le\hat{I}_3\ \text{ for all }r\in[0,t]\big\}.
	\]
	$\mathcal J$ is an interval containing $0$, since
	$|I_3^0|\le\hat{I}_3$, and it is closed in $[0,\omega)$ by continuity
	of $I_3$.  Fix $t\in\mathcal J$.  On $[0,t]$ we have
	$|I_3|\le\hat{I}_3$ and, by Step 1, $|I_2|\le\hat{I}_2$.  With these
	bounds in place of Assumption~\ref{ass:compat} and with
	\eqref{eq:plantb}, Steps 2--4 and~6 of the proof of
	Theorem~\ref{thm:main} apply with $[0,t]$ in place of $[0,\omega)$,
	the norms of $k_1$, $\dot{e}_1$ and $B$ being those of
	$L^\infty([0,t];\mathbb{R})$.  They give
	$\|k_1\|_\infty\le\bar{k}_1$ by \eqref{eq:e1bound} and
	$\|\dot{e}_1\|_\infty\le\hat{P}_1$ by \eqref{eq:pbound},
	\eqref{eq:e1}, \eqref{eq:I1} and $|I_1|\le\hat{I}_1$; hence, by
	\eqref{eq:MB}, \eqref{eq:v1bar}, \eqref{eq:chain} and \eqref{eq:ap1},
	\[
	|v_1(r)|\le\max\big\{|v_1(0)|,\|B\|_\infty\big\}\le\bar{v}_1
	<\underline{M}_1+\kappa_3\hat{I}_3/\nu_1
	\qquad\text{for all }r\in[0,t],
	\]
	which for $t=0$ is immediate from \eqref{eq:chain} and
	\eqref{eq:ap1}.  By continuity of $v_1$ there is $\eta>0$ such that
	$|v_1|<\underline{M}_1+\kappa_3\hat{I}_3/\nu_1$ on
	$[0,t+\eta)\cap[0,\omega)$, hence
	$\nu_1|w_1|\le\nu_1(|v_1|-\underline{M}_1)_+\le\kappa_3\hat{I}_3$ there,
	and variation of constants in \eqref{eq:I2} gives
	\[
	|I_3(r)|\le\mathrm e^{-\kappa_3r}|I_3^0|
	+\big(1-\mathrm e^{-\kappa_3r}\big)\hat{I}_3\le\hat{I}_3
	\qquad\text{for all }r\in[0,t+\eta)\cap[0,\omega).
	\]
	Hence $\mathcal J$ is also open in $[0,\omega)$, thus
	$\mathcal J=[0,\omega)$.

	Both bounds hold along every solution, in particular along every
	maximal one, which is Assumption~\ref{ass:compat}.
\end{proof}

\begin{remark}\label{rem:apread}
	Condition \eqref{eq:ap2} is the two-sided analogue of the
	saturation-level bound of \cite{HopIlcRya10}.  The rate limit does
	not enter it; it enters the constant \eqref{eq:ap2bound} through
	$\bar{W}_2$, the actuator needing the time $W_2/\hat{r}$ to traverse
	the constraint set. In
	\eqref{eq:ap1}, $\bar{v}_1$ depends on $\hat{I}_3$ by
	\eqref{eq:chain}; for $\hat{I}_3=0$ the condition reads
	$\bar{v}_1<\underline{M}_1$, so that $u_1$ does not saturate, and for
	$\hat{I}_3>0$ saturation is admitted.
\end{remark}

%=====================================================================
\section{Simulation study}\label{sec:sim}
%=====================================================================

Three questions are studied on the nonlinear model: how the controller
performs on the ITER benchmark, what the input constraints cost when
they become active, and how tight the compatibility conditions are.

\subsection{Setup}\label{sec:setup}

We integrate the nonlinear model \eqref{eq:model}, \eqref{eq:closure} at
the operating point of Section~\ref{sec:closure} with
$T_c=1\,\mathrm{s}$ over a flat-top of $400\,\mathrm{s}$, by the
classical Runge--Kutta method with fixed step $4\,\mathrm{ms}$
($1\,\mathrm{ms}$ in Section~\ref{sec:apcheck}, where every pulse edge
falls on the grid) and without event detection, the discontinuous $d_4$
being sampled at the Runge--Kutta stages. Halving the step three
times changes the benchmark errors by $0.001$ percentage points.
The initial state is perturbed by $+1\%$ in $n_{DT}$ and $-0.5\%$ in
$E$; the controller starts from $\xi(0)=I_2(0)=I_1(0)=I_3(0)=0$ and
$u_2(0)=0$.  All disturbances of
Section~\ref{sec:dist} are active, with the indicator function
$\mathbf 1$ and a factor $m>0$ scaling the particle sources, $m=1$ on
the benchmark,
\begin{equation}\label{eq:distspec}
	\begin{aligned}
		d_1(t)&=1+0.015\sin\tfrac{2\pi t}{140\,\mathrm s}
		-0.03\,\mathbf 1_{[260,290)}(t),\\
		d_2(t)&=0.010\,m\,\tfrac{\bar n_\alpha}{\bar\tau}
		\big(1+\tfrac12\sin\tfrac{2\pi t}{95\,\mathrm s}\big),\qquad
		d_3(t)=1.2\,d_2(t)\,\bar n_{DT}/\bar n_\alpha,\\
		d_4(t)&=-\tfrac{3\cdot10^{-4}\bar E}{50\,\mathrm{ms}}
		\mathbf 1_{[0,50\,\mathrm{ms})}(t\bmod2\,\mathrm s)
		-\tfrac{4\cdot10^{-4}\bar E}{100\,\mathrm{ms}}
		\mathbf 1_{[0,100\,\mathrm{ms})}(t\bmod20\,\mathrm s).
	\end{aligned}
\end{equation}
For $m=1$ the particle sources $d_2,d_3$ are at $7\%$ and $3.6\%$ of
the nominal transport losses, and $d_4$ removes $0.03\%$ of the nominal
stored energy $\bar E$ every $2\,\mathrm s$ and $0.04\%$ every
$20\,\mathrm s$.  These bursts presuppose mitigated events: ITER relies on ELM
control to reduce the loss per ELM from about $20\,\mathrm{MJ}$, i.e.\ $6\%$ of
$V\bar E=332\,\mathrm{MJ}$, to about $0.7\,\mathrm{MJ}$ \cite{Loa10}, and a single burst of a few per cent
would exceed the asymptotic width $2.5\%$ of the funnel chosen below; bursts of
$0.3\%$ and $0.4\%$ are treated in Section~\ref{sec:apcheck}.  Since $d_2,d_3\ge0$, the standing
hypothesis of Proposition~\ref{prop:membership}\,(ii) holds.
Figures in MW are totals, $V\bar P=36.4\,\mathrm{MW}$.  The
simulation works in the inputs of Section~\ref{sec:inputconstr} and
the normalised outputs $y_1=n_e/\bar n_e$, $y_2=E/\bar E$, a rescaling
of \eqref{eq:class} under which
$\gamma=\bar S/(\bar n_ek_d\tau_E)$, $a_2=\bar P/\bar E
=0.1096\,\mathrm{s^{-1}}$ and $\bar B$ is unchanged; every quantity
of Sections~\ref{sec:main} and \ref{sec:apriori} quoted below is
evaluated in these units.

The constraint sets of \eqref{eq:sat} are
$P_{\mathrm{aux}}(t)\in[0,73]\,\mathrm{MW}$, according to the installed capacity
\cite{Shi07}, and $S(t)\in[0,3\bar S]$, i.e.\ $\underline u_i=-1$,
$\overline u_1=2$, $\overline u_2=1.004$ and
$\underline M_1=\underline M_2=1$; the upper fuelling limit is never
reached.  The heating power is rate limited with
$\bar r=10\,\mathrm{MW\,s^{-1}}$; without \eqref{eq:act} the law
would demand up to $32\,\mathrm{MW\,s^{-1}}$ on the benchmark.  The controller \eqref{eq:controller} uses the references
$y_{\mathrm{ref},1}=y_{\mathrm{ref},2}=1$,
$\psi_1(t)=0.10\mathrm{e}^{-t/60\,\mathrm{s}}+0.030$,
$\psi_2(t)=0.10\mathrm{e}^{-t/60\,\mathrm{s}}+0.025$ on the relative
errors, $\hat\theta=0.10$, $\kappa_1=\kappa_2=\kappa_3=1$,
$\nu_1=\nu_2=0.1$ and $k_3=20$.   The barrier arguments
$\varphi_i|e_i|$ and $|\theta|/\hat\theta$ are clipped at $0.999$; from
the instant the clip is active the integrated system is a
\emph{regularised} controller and no longer \eqref{eq:controller}.  We
report that instant, which no run of \eqref{eq:controller} below reaches;
runs approaching the clip are stiff and use steps of at most
$0.025\,\mathrm{ms}$.

\begin{figure}[t]
	\begin{center}
		\includegraphics[width=0.92\textwidth]{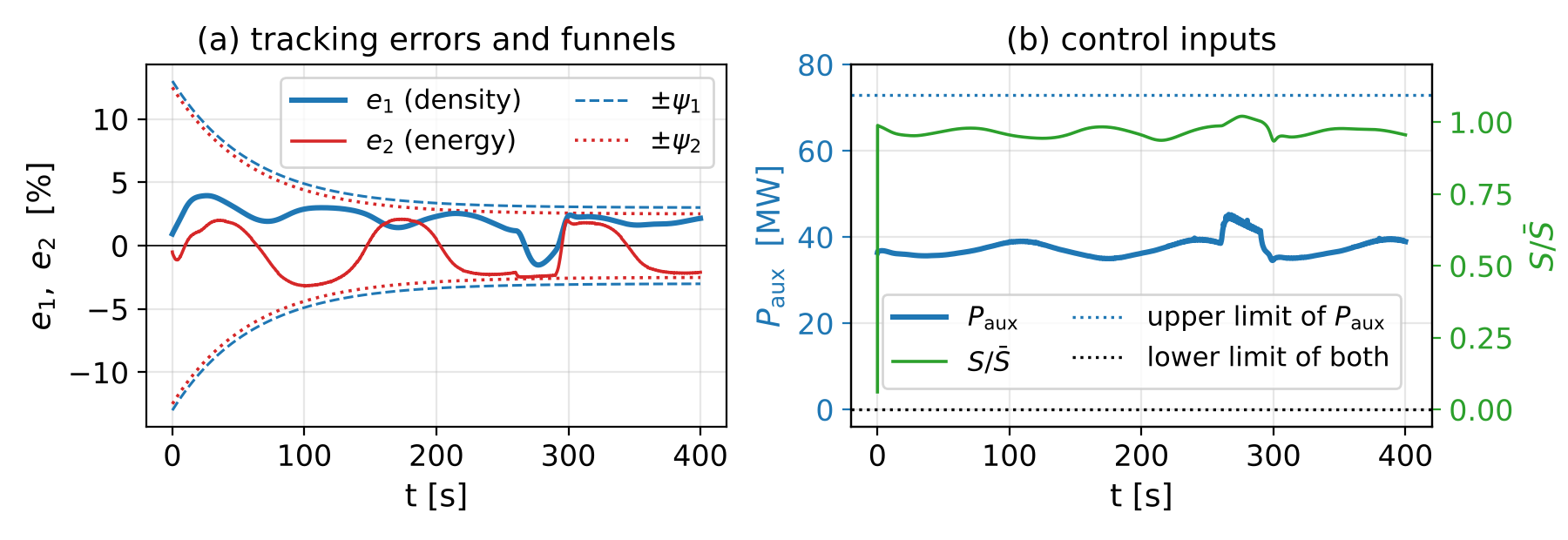}
		\vspace{-4mm}
		\caption{ITER benchmark with the controller \eqref{eq:controller}.
			(a) The errors $e_1$ and $e_2$ with the funnel boundaries
			$\pm\psi_1$ and $\pm\psi_2$; on the benchmark $I_1=I_2=0$, so that
			$e_1=n_e/\bar n_e-1$ and $e_2=E/\bar E-1$.  (b) Auxiliary heating
			power (left axis) and fuelling rate (right axis) with their limits
			(dotted).}
		\label{fig:bench}
	\end{center}
\end{figure}

\subsection{Results}

Figure~\ref{fig:bench} shows the outcome.  The
proposed controller keeps all three errors inside their funnels, and
neither amplitude limit is reached, so $I_1,I_2,I_3$ vanish identically
and Theorem~\ref{thm:main}\,(iii) gives the prescribed funnels.  The
deviations $\delta n_e(t):=n_e(t)-\bar n_e$ and
$\delta E(t):=E(t)-\bar E$ from the operating point satisfy
$\max_t|\delta n_e(t)|/\bar n_e=3.9\%$ and
$\max_t|\delta E(t)|/\bar E=3.2\%$, while $P_{\mathrm{aux}}$ ranges
between $34.6$ and $45.2\,\mathrm{MW}$ and $S/\bar S$ between $0.06$
and $1.02$.  The
rate limit is reached in short spikes at the energy bursts $d_4$, at a
cost of under $0.01$ percentage points in either error.

A model-based state feedback as in Section~\ref{sec:related} requires
all five states, among them the unmeasured $n_n$ and
$\tau_E$, and the model \eqref{eq:model}, \eqref{eq:closure}, and provides no
prescribed bounds on the transient; the proposed law uses the two diagnosed
outputs, no parameter, and returns the funnels of
Fig.~\ref{fig:bench}(a).

\subsection{When the actuators reach their limits}\label{sec:influx}

The benchmark reaches no amplitude limit.  We therefore increase the
factor $m$ in \eqref{eq:distspec},
which drives the constraint $S\ge0$: an influx excess must be answered
by reducing the fuelling, which cannot go below zero. We further compare the controller \eqref{eq:controller} with its variant (P), in which the corrections are
switched off, $I_2\equiv I_1\equiv I_3\equiv0$, the demanded input is
clipped at the actuator limits and \eqref{eq:act} is dropped.
While no amplitude limit is reached, (P) is \eqref{eq:controller}
without \eqref{eq:act}. Table~\ref{tab:sat} reports the outcome for \eqref{eq:controller} and (P).

\begin{table}[t]
	\caption{Controller \eqref{eq:controller} with the particle influx scaled by $m$, all other
		disturbances and both funnels as in Section~\ref{sec:setup}:
		$\max|\delta n_e/\bar n_e|$, the bound
		$\max_t\big(\psi_1(t)+|I_1(t)|\big)$ from
		Theorem~\ref{thm:main}\,(iii) and the observed
		$\hat{I}_3^{\,\mathrm{obs}}=\max_t|I_3(t)|$, from which $\max|I_1|$
		differs by less than $0.005$; $e_1$ never leaves its funnel.  Last row: first instant at which the
		clip of Section~\ref{sec:setup} becomes active in (P).}
	\label{tab:sat}
	\begin{center}\small
		\begin{tabular}{@{}lrrrrrrrr@{}}
			\hline
			$m$ & $1$ & $6$ & $10$ & $12$ & $12.5$ & $15$ & $18$ & $20$\\
			\hline
			$\max|\delta n_e/\bar n_e|$ [\%] & $3.9$ & $10.2$ & $11.3$ & $12.2$ & $12.4$ & $14.2$ & $28.0$ & $46.6$\\
			bound [\%] & $13.0$ & $13.0$ & $13.0$ & $13.0$ & $13.2$ & $14.8$ & $28.2$ & $46.7$\\
			$\hat{I}_3^{\,\mathrm{obs}}$ & $0$ & $0$ & $0.001$ & $0.010$ & $0.013$ & $0.029$ & $0.195$ & $0.383$\\
			clip of (P) [s] & -- & -- & -- & -- & $4.1$ & $2.9$ & $2.3$ & $2.0$\\
			\hline
	\end{tabular}\end{center}
\end{table}

\begin{figure}[t]
	\begin{center}
		\includegraphics[width=0.85\textwidth]{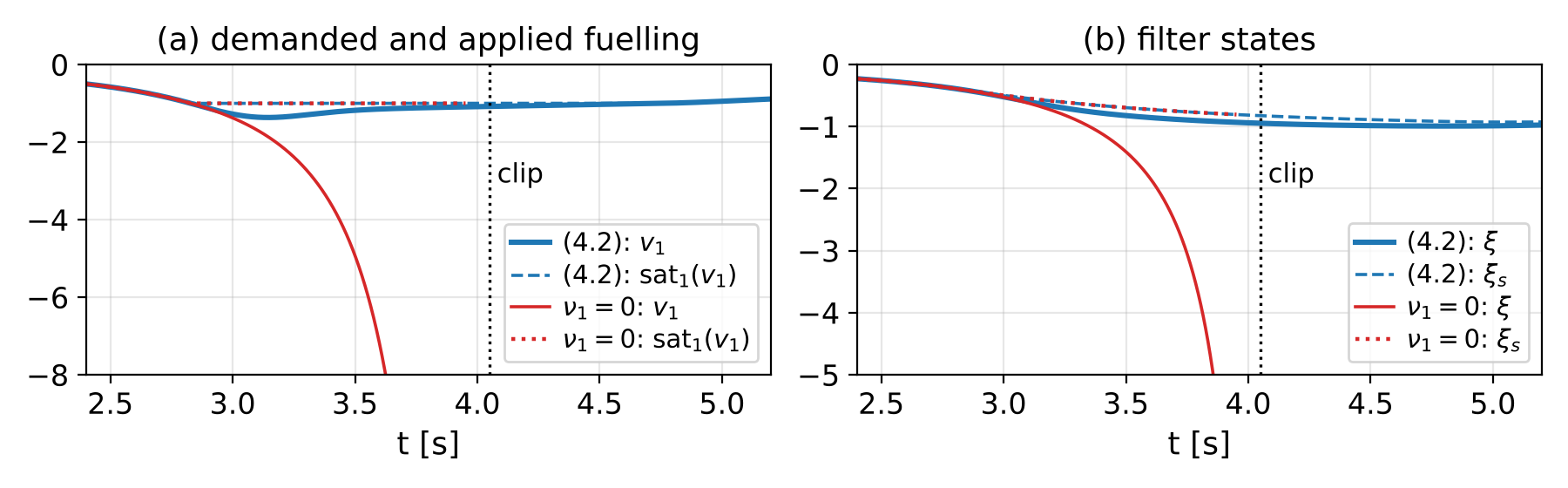}
		\vspace{-4mm}
		\caption{First saturation episode for $m=12.5$: controller
			\eqref{eq:controller} and its variant with $\nu_1=0$.  Here $\xi_s$ is the state
			of the filter \eqref{eq:xi} driven by the applied input,
			$\dot\xi_s=-\kappa_3\xi_s+\operatorname{sat}_1(v_1)$,
			$\xi_s(0)=\xi^0$.  Dotted line: first clip of the
			variant, whose curves are drawn up to $t=3.95\,\mathrm s$.}
		\label{fig:mech}
	\end{center}
\end{figure}

Three things are visible.  First, (P) reaches the clip for
$m\ge12.5$, and the density correction alone decides this:
Figure~\ref{fig:mech} shows, for $m=12.5$, the variant of
\eqref{eq:controller} with $\nu_1=0$, in which $I_1\equiv I_3\equiv0$
and the heating channel is unchanged.  From $t=2.83\,\mathrm s$ the
fuelling sits at its lower limit.  Without the correction the demand
$v_1$ and the filter state $\xi$, which \eqref{eq:xi} drives by the
demand, depart from $\operatorname{sat}_1(v_1)$ and $\xi_s$, and the
clip is reached at $t=4.05\,\mathrm s$, as in (P).  With the
correction, $v_1$ stays within $0.37$ of the limit and
$|\xi-\xi_s|=|I_3|/\nu_1\le0.14$ (Lemma~\ref{lem:shadow}), and the
episode ends at $t=4.75\,\mathrm s$.  At $t=4.05\,\mathrm s$ the
density deviation is $12.3\%$ in both runs: the correction does not
reduce the excursion, it keeps the law defined.  Beyond the clip the
regularised system is stiff and its data depend on the integration
step, so none are reported.

Second, \eqref{eq:controller} keeps the error inside its funnel for
every $m$.  At $m=20$, $\max|\theta|=0.40\hat\theta$, $\max|v_1|=4.83$
and $\max\varphi_1|e_1|=0.996$.

Third, the bound $\psi_1+|I_1|$ of Theorem~\ref{thm:main}\,(iii) is
within $0.2$ percentage points of the realised excursion for $m\ge18$.
The operating point sits at $85\%$ of the Greenwald limit, $18\%$ of
headroom; the observed band stays inside that headroom for $m\le15$ and
exceeds it for $m\ge18$.  The a priori band $\psi_1+\hat{I}_1$ of
Theorem~\ref{thm:main}\,(iii) needs the constant $\hat{I}_3$, for which
these runs supply only the observed value.

\subsection{Checking the compatibility conditions}\label{sec:apcheck}

We evaluate the compatibility conditions on the benchmark and on a
\emph{raised} disturbance set -- $0.3\%$ ELM loss at $1\,\mathrm{Hz}$, $0.4\%$
sawtooth loss every $20\,\mathrm s$, a $5\%$ confinement modulation
with a $10\%$ drop for $40\,\mathrm s$, particle influx tripled.  There
the proposed law keeps both funnels with $8.6\%$ and $5.7\%$ maximum
error and $\hat{I}_3^{\,\mathrm{obs}}=0.001$; (P) reaches the clip in
the energy channel at $t=240\,\mathrm s$.  We evaluate condition
\eqref{eq:ap2} of Proposition~\ref{prop:ap} for the energy channel on
$[0,T_f]$, $T_f=400\,\mathrm s$, in which
$\underline M_2=|\underline u_2|=1$, $\overline u_2=1.004$ and the authority
is $a_2\underline M_2=1.10\cdot10^{-1}\,\mathrm{s^{-1}}$; the references being
constant, the drift is $\max_t|F_2(t)|$.  Except for $\hat I_2$, all quantities below are
\emph{observed}, none an a priori bound.

On the benchmark the drift is $2.99\cdot10^{-2}\,\mathrm{s^{-1}}$, below the
authority by a factor of $3.7$, and the demand never leaves the box, so
$I_2\equiv0$.  The constant
$\hat I_2$ of \eqref{eq:ap2bound}, which does not depend on the disturbances, is nevertheless $116\%$: its first
term is $81\%$, the actuator needing $W_2/\hat r=7.3\,\mathrm s$ to
traverse the constraint set, the second is $\|\psi_2\|_\infty=12.5\%$,
and the third is $22\%$.  The plant belongs to the class only for
$\hat I_2<87.5\%$, so Proposition~\ref{prop:ap} does not apply to the
benchmark.

Under the raised disturbances the drift reaches
$1.87\cdot10^{-1}\,\mathrm{s^{-1}}$ and exceeds the authority for
$1.0\%$ of $T_f$, by a factor of up to $1.71$, so \eqref{eq:ap2} fails.
The demand is outside its box for $6.3\%$ of $T_f$ ($5.7\%$ above
$\overline u_2$, $0.5\%$ below $\underline u_2$), in $165$ episodes, the
longest of $5.8\,\mathrm s$; the observed shift stays at
$\max|I_2|=1.3\%$.  %A bound for a drift that exceeds the authority only in pulses is not provided by Section~\ref{sec:apriori}.

Condition \eqref{eq:ap1} holds for no $\hat{I}_3$ on the benchmark:
already with $M_1^F=0$, \eqref{eq:chain} gives $\bar{k}_1\ge71$ and
$\bar{v}_1\ge5\cdot10^{4}$ for every $\hat{I}_3\ge0$, so that
\eqref{eq:ap1} asks for $\hat{I}_3>5\cdot10^{3}$, against the range
$\hat{I}_3<0.87$ in which the plant belongs to the class and an observed
$\max|v_1|=0.94<\underline{M}_1=1$ -- the high-gain conservatism
noted in \cite[Rem.~7]{FotRov24}.

%=====================================================================
\section{Conclusions}\label{sec:conclusions}
%=====================================================================

For the zero-dimensional burn model with a dynamic, loss-power based
confinement closure, the two-actuator configuration has vector relative degree $(2,1)$ (Proposition~\ref{prop:struct}) and admits
a funnel controller that needs no output derivative and meets every actuator constraint by
construction.  The invariance of the funnels (Theorem~\ref{thm:main})
needs no sign condition on the input gain, but it is conditional: the plant re-enters through
the compatibility assumption, which asks that the reference corrections
stay bounded. A sufficient (yet conservative) a priori condition for the compatibility assumption is derived in Proposition~\ref{prop:ap}. We note that profile effects \cite{GraSch25}, measurement errors and
pulsed particle losses, which (D2) excludes, are outside the present
scope and potential topics for future research.

\section*{Acknowledgements}

Anthropic Claude was used to assist with language editing of the manuscript and with the development of parts of the numerical implementation. All mathematical derivations, numerical results, and scientific conclusions are those of the author. The author independently verified all AI-assisted code and reviewed all AI-assisted content, and remains fully responsible for the content of the manuscript.

\end{document}